\documentclass[onefignum,onetabnum]{siamart251216}
\usepackage[normalem]{ulem}
\newif\ifshowchanges \showchangestrue

\usepackage{amssymb}
\usepackage{mathtools}
\usepackage{bm}
\usepackage{booktabs}
\usepackage{array}
\usepackage{enumitem}
\usepackage{algorithmic}
\usepackage{epstopdf}
\ifpdf
  \DeclareGraphicsExtensions{.eps,.pdf,.png,.jpg}
\else
  \DeclareGraphicsExtensions{.eps}
\fi

\newsiamthm{assumption}{Assumption}
\newsiamremark{remark}{Remark}

\newcommand{\R}{\mathbb R}
\newcommand{\E}{\mathbb E}
\newcommand{\V}{\mathbb V}
\newcommand{\Hh}{\mathbb H}
\newcommand{\C}{\mathcal C}
\newcommand{\I}{\mathcal I}

\newcommand{\Mcal}{\mathcal M}
\newcommand{\Kcal}{\mathcal K}

\newcommand{\norm}[1]{\left\lVert #1\right\rVert}
\newcommand{\Proj}{\operatorname{Proj}}
\newcommand{\dom}{\operatorname{dom}}

\usepackage[normalem]{ulem}
\newcommand{\revmode}{1}
\DeclareRobustCommand{\rev}[1]{%
  \ifcase\revmode\relax\or{\color{red}\relax#1}\else#1\fi}
\DeclareRobustCommand{\revdel}[1]{%
  \ifcase\revmode\relax#1\or{\color{red}\relax\sout{#1}}\kern0.15em\else\fi}
\DeclareRobustCommand{\revdelm}[1]{%
  \ifcase\revmode\relax#1\or\text{\color{red}\sout{$#1$}}\else\fi}
\DeclareRobustCommand{\revid}[1]{%
  \ifnum\revmode=1\relax\textsuperscript{\color{red}\scriptsize[#1]}\fi}
\headers{Box-Constrained Elliptic Control under Uncertainty}
        {Z. Li, A. Onwunta, and H. Wang}

\title{A Preconditioned Semismooth Newton Method for  Box-Constrained Elliptic Optimal Control Problems with Random Coefficients
}

\author{Zhendong Li\thanks{Department of Industrial and Systems
  Engineering, Lehigh University, Bethlehem, PA
  (\email{zhl923@lehigh.edu}, \email{ako221@lehigh.edu}).}
\and Akwum Onwunta\footnotemark[1]
\and Hailing Wang\thanks{Department of Applied Mathematics, The Hong Kong Polytechnic University, Hung Hom, Kowloon, Hong Kong (\email{hailing.wang@polyu.edu.hk}).}}

\ifpdf
\hypersetup{
  pdftitle={A Preconditioned Semismooth Newton Method for  Box-Constrained Elliptic Optimal Control Problem with Random Coefficients},
  pdfauthor={Z. Li, A. Onwunta, and H. Wang}
}
\fi

\begin{document}

\maketitle

\begin{abstract}
We consider the optimal control problem of an elliptic PDE with random diffusion under pointwise box constraints on the control. The constraints make the optimality conditions nonsmooth, with an active set that is not known in advance, and the random diffusion couples the state and the adjoint across all stochastic modes. To design a solution algorithm for this problem, we rely on the Fenchel duality viewpoint.  Fenchel duality absorbs the box constraint into the objective and yields an exact reduced problem in the elliptic adjoint variable alone.  After a stochastic Galerkin--finite element discretization, eliminating the state leaves a strongly monotone, piecewise smooth equation, which we solve directly by a matrix-free semismooth Newton method with the box projection retained exactly and with no outer iteration, penalty or smoothing parameter. Next, we introduce and analyze a current-set preconditioner, which adds the active-set correction of the current Newton step to a mean-diffusion reference operator and is inverted exactly by a rank-1 update. For implementation, we use a cheaper preconditioner that keeps only a fraction of free control values. Numerical experiments  show second-order convergence of the state, and Newton and PCG counts that stay bounded under mesh refinement and growth of the stochastic dimension.
\end{abstract}



\section{Introduction}\label{sec:intro}
 Optimal control problems governed by partial differential equations (PDEs) are ubiquitous in computational science and engineering, and both a mature theory and a
large body of numerical methods are available in the deterministic
case~\cite{Troltzsch-2010-OCP,elman_finite_2014}. Pointwise control bounds $a\le u\le b$ change the picture.
The first-order conditions are no longer a system
of equations but a variational inequality, or equivalently a nonsmooth equation
built from the metric projection onto the admissible
box~\cite{Ito-2008-LMA,Facchinei-2003-FDV}.  Three consequences follow.
First, although the reduced objective is strongly convex, the resulting
Karush--Kuhn--Tucker (KKT) system is nondifferentiable, so classical Newton
methods do not apply.  Second, the set of indices at which the bounds are
attained---the active set---is not known in advance and must be discovered as
part of the solution process. Third, the standard remedies---penalization, interior-point
regularization, smoothing---replace the exact constraint by an approximation
governed by an additional parameter, and recover it only in the limit.

In the deterministic problem,
the primal--dual active set (PDAS) strategy of Bergounioux, Ito and
Kunisch~\cite{Bergounioux-1999-PDS} treats the bounds exactly by alternating
between a prediction of the active set and the solution of an
equality-constrained subproblem.  Hinterm\"uller, Ito and
Kunisch~\cite{Hintermuller-2002-PDA} subsequently showed that PDAS is 
precisely a semismooth Newton (SSN) method for the projection-based reformulation, which explains its local superlinear and mesh-independent convergence~\cite{Clarke-1983-ONA,Mifflin-1977-SFN,Qi-1993-NNM,Ulbrich-2011-SNM,Hintermuller-2004-MIS}.

In most applications the PDE coefficients are uncertain and are modelled by random fields and lead to optimal control problems under uncertainty~\cite{Martinez-Frutos-2018-OCP}. The random variable may be discretized by sampling, by stochastic collocation~\cite{Chen-2014-WRB,Kouri-2013-MSC,Kouri-2013-TRA,Tiesler-2012-SCO}, or by stochastic Galerkin (SG) projection onto a global orthogonal polynomial basis~\cite{Ghanem-1991-SFE,Xiu-2010-NMS,LeMaitre-2010-SMU,Lord-2014-ICS,%
Babuska-2004-GFE,Lee-2013-SGM,Rosseel-2012-OCS}. We use SG, which yields a single coupled system: the discrete state operator is a sum of Kronecker products, and state and adjoint interact across all stochastic modes.

Combining pointwise constraints with random coefficients is harder than
treating either feature alone.  
For {control}
constraints, the projection onto the box is available in closed form and need
not be regularized, yet stochastic Galerkin discretizations of such problems
have been studied mainly from the approximation side: Shen et
al.~\cite{Shen-2015-SGC} derive a priori error estimates for obstacle-type
control constraints, and Winkler and Y\"ucel~\cite{Winkler-2025-SGD} combine
error estimates for a Dirichlet boundary control problem with a primal--dual
active set solver and block-diagonal preconditioning. There, the control is
deterministic and acts on $\partial D$, so the active set stays decoupled from
the stochastic tensor structure. Solvers for stochastic control problems are developed in~\cite{Benner-2016-BDP, BDOS2016, HDO2023,li2025stochastic, Antil2025, Antil-2025-SMY, BENNER2020}. None of them treats control constraints.

An alternative is the  Fenchel duality viewpoint, where the dual absorbs the constraint into the objective.  
Clason and Kunisch~\cite{Clason-2011-DBA} found the duals of certain unconstrained
elliptic control problems structurally more favourable than the primal
problems, and Wang et al.~\cite{MR4750911} derived the Fenchel dual of a
linear parabolic problem with box constraints on the control, in which the
constraints are embedded in an unconstrained, differentiable dual objective,
and solved it by conjugate gradients and by semismooth Newton.  They also
identified the obstacle that reappears here: semismooth Newton applied to the
primal problem treats the box and the PDE constraint simultaneously, so the
varying active set changes the structure of the Newton equation and the
preconditioner has to be adjusted at every step. 
In the stochastic setting, this obstacle is sharper: control values are indexed by pairs of a spatial node and a quadrature point, so the active-set term of the Jacobian matrix has no Kronecker structure. A preconditioner containing it exactly must be rebuilt at every step, whereas a mean-based Kronecker preconditioner is cheap but blind to the location of the free set.

In this paper, we develop a dual-based SSN framework for
optimal control problems governed by elliptic PDEs with random coefficients
and pointwise box constraints on the control.
We derive the Fenchel dual in Bochner spaces over the parameter domain,
which absorbs the control constraints into an unconstrained dual objective
in the adjoint variable, and establish the corresponding optimality system
and primal--dual recovery relations.
After stochastic Galerkin--finite element discretization, we derive the
exact finite-dimensional dual and solve the resulting reduced equation
directly by a matrix-free semismooth Newton--PCG method, in which the box
projection is retained exactly, and no outer iteration, penalty, or smoothing parameter is required.
To handle the active set, which changes at every Newton step, we analyze a current-set preconditioner $P_W$ that places the
current active-set correction beside a mean-based
reference operator and applies its inverse exactly through the
Sherman--Morrison--Woodbury (SMW) identity~\cite{Hager-1989-UIM}.
For numerical implementation, we use a cheap approximation  $P_k$ of the preconditioner $P_W$, which keeps only the
fraction of free control values, and its application needs one sparse
Cholesky factor of a spatial matrix, which is reused across Newton steps.

{The remainder of the paper is organized as follows.
{Section~\ref{sec:problem}-\ref{sec:fenchel} presents the problem and derives
the Fenchel dual and the optimality system.
Section~\ref{sec:sgfe} gives the SG-FE discretization and its exact dual.
Section~\ref{sec:ssn} develops the direct SSN--PCG method and
the spectral analysis of the current-set preconditioner $P_W$.
Section~\ref{sec:numerics} specifies the practical preconditioner $P_k$  and reports the numerical experiments.
Section~\ref{sec:conclusion} draws conclusions and discusses future work.} 

\section{Problem formulation and well-posedness}\label{sec:problem}

\subsection{The optimal control problem}

Let \(D\subset\R^d\), \(d\in\{1,2,3\}\), be a bounded Lipschitz domain with
Poincar\'e constant \(C_P\), that is,
\(\norm{v}_{L^2(D)}\le C_P\norm{\nabla v}_{L^2(D)^d}\) for all
\(v\in H_0^1(D)\), let $\xi=(\xi_1,\dots,\xi_K)\colon(\Omega,\mathcal F,\mathbb P)\to\R^K$, $K<\infty$, be a random vector with image $\Gamma_\xi\coloneqq\xi(\Omega)$, and let \((\Gamma_\xi,\mathcal B,\rho)\) be the
induced parameter probability space, where $\mathcal B$ is the Borel $\sigma$-algebra on~$\Gamma_\xi$ and $\rho\coloneqq\mathbb P\circ\xi^{-1}$ is the law of~$\xi$.
Let~\(y\) denote the state and \(u\)  the
control, and set
\[
 V=H_0^1(D),\qquad H=L^2(D),\qquad
 \V=L^2_\rho(\Gamma_\xi;V),\qquad
 \Hh=L^2_\rho(\Gamma_\xi;H),
\]
where the norms and inner products on \(\V\) and \(\Hh\) include the
expectation with respect to \(\rho\). We equip \(V\) with the norm
\(\norm{v}_V\coloneqq\norm{\nabla v}_{L^2(D)^d}\), so that
\(\norm{v}_{\V}^2=\int_{\Gamma_\xi}\norm{\nabla v(\xi)}_{L^2(D)^d}^2\,d\rho(\xi)\)
and \(\norm{v}_{\Hh}\le C_P\norm{v}_{\V}\) for \(v\in\V\). Given a desired state \(y_d\in\Hh\), a
regularization parameter \(\gamma>0\), finite bounds \(a<b\), and a random
diffusion coefficient \(\kappa\), the problem studied in this paper is
\begin{equation}
\label{eq:ocp}
\begin{aligned}
 \min_{y,\,u}\quad
 &\mathrlap{\frac12\int_{\Gamma_\xi}\!\int_D\bigl|y(x,\xi)-y_d(x,\xi)\bigr|^2dx\,d\rho(\xi)
+\frac\gamma2\int_{\Gamma_\xi}\!\int_D\bigl|u(x,\xi)\bigr|^2dx\,d\rho(\xi)}\\
 \text{subject to}\quad
 &-\nabla\cdot\bigl(\kappa(x,\xi)\nabla y(x,\xi)\bigr)=u(x,\xi)
   &&\text{in }D\times\Gamma_\xi,\\
 &y(x,\xi)=0
   &&\text{on }\partial D\times\Gamma_\xi,\\
 &a\le u(x,\xi)\le b
   &&(dx\otimes d\rho)\text{-a.e.}
\end{aligned}
\end{equation}

Both integrals are taken with respect to the product of the Lebesgue measure
on \(D\) and the probability measure \(\rho\); so, the objective in
\eqref{eq:ocp} is exactly
\[
 \frac12\norm{y-y_d}_{\Hh}^2+\frac\gamma2\norm{u}_{\Hh}^2,
\]
and we use this compact form from Section~\ref{sec:fenchel} onwards.

We  note here that 
the
bounds are imposed pointwise in the spatial and the parameter variables,
so the admissible set
\begin{equation}
\label{eq:admissible}
 \C=\left\{u\in\Hh:
 a\le u(x,\xi)\le b
 \text{ for }(dx\otimes d\rho)\text{-a.e. }(x,\xi)\right\}
\end{equation}
is a subset of the full Bochner space rather than of a deterministic space of
controls.  Since \(D\) and the parameter space carry finite measure, \(\C\) is
nonempty, bounded, closed, and convex in \(\Hh\).  


\subsection{The state operator and the reduced problem}

Assume that \(\kappa\) is jointly measurable, symmetric, and uniformly
elliptic:

\begin{equation}
\label{eq:UE}
 0<\kappa_-\le\kappa(x,\xi)\le\kappa_+<\infty,
\end{equation}
for all \(\eta\in\R^d\) and almost every \((x,\xi)\), with deterministic
constants \(\kappa_-\) and \(\kappa_+\).  The weak form of the
state equation in \eqref{eq:ocp} is governed by the bilinear form
\[
 \mathfrak a(y,v)=
 \int_{\Gamma_\xi}\int_D
 \kappa(x,\xi)\nabla y(x,\xi)\cdot\nabla v(x,\xi)
 \,dx\,d\rho(\xi).
\]

\begin{theorem}[Lax--Milgram {\cite{Troltzsch-2010-OCP}}]
\label{thm:lax-milgram}
Let \(W\) be a real Hilbert space and let \({\mathfrak t}:W\times W\to\R\) be a
bilinear form for which there are constants \(M\ge\alpha>0\) such that
\[
 |{\mathfrak t}(w,z)|\le M\norm{w}_W\norm{z}_W
 \quad\text{and}\quad
 {\mathfrak t}(w,w)\ge\alpha\norm{w}_W^2
 \qquad(w,z\in W).
\]
Then, for every \(\ell\in W'\), there is a unique \(z\in W\) such that
\[
 {\mathfrak t}(z,w)=\ell(w)\qquad(w\in W),
\]
and \(\norm{z}_W\le\alpha^{-1}\norm{\ell}_{W'}\).
\end{theorem}

\begin{proposition}
\label{prop:state-operator}
Under \eqref{eq:UE}, for every \(u\in\Hh\) there is a unique \(y=Su\in\V\)
such that
\begin{equation}
\label{eq:state-weak}
 \mathfrak a(y,v)=(u,v)_{\Hh}\qquad(v\in\V).
\end{equation}
The solution map
\[
 S:\Hh\longrightarrow\V\hookrightarrow\Hh
\]
is bounded, linear, positive, injective, and self-adjoint as an operator on
\(\Hh\).
\end{proposition}

\begin{proof}
Boundedness and coercivity of \(\mathfrak a\) on \(\V\) follow from
\eqref{eq:UE}.  Moreover,
\[
 |(u,v)_{\Hh}|
 \le\norm{u}_{\Hh}\norm{v}_{\Hh}
 \le C_P\norm{u}_{\Hh}\norm{v}_{\V}.
\]
Theorem~\ref{thm:lax-milgram} therefore gives a unique solution and
\[
 \norm{Su}_{\V}\le\frac{C_P}{\kappa_-}\norm{u}_{\Hh},
 \qquad
 \norm{Su}_{\Hh}\le\frac{C_P^2}{\kappa_-}\norm{u}_{\Hh}.
\]
Linearity of \(S\) follows from uniqueness. 
For \(u_1,u_2\in\Hh\) and \(c_1,c_2\in\R\), the bilinearity of \(\mathfrak a\) and of \((\cdot,\cdot)_{\Hh}\) shows that \(c_1 Su_1+c_2 Su_2\) solves~\eqref{eq:state-weak} with right-hand side \(c_1 u_1+c_2 u_2\). By uniqueness, \(S(c_1 u_1+c_2 u_2)=c_1 Su_1+c_2 Su_2\), that is, \(S\) is linear.
For \(f,g\in\Hh\), utilizing the symmetry of this form, we obtain from the equation corresponding to \(Sg\) (with \(Sf\)) and the equation corresponding to \(Sf\) (with \(Sg\)) that
\[
 (Sf,g)_{\Hh}=(g,Sf)_{\Hh}=\mathfrak a(Sg,Sf)=\mathfrak a(Sf,Sg)=(f,Sg)_{\Hh}.
\]
Hence \(S^*=S\).  Moreover, \((Sf,f)_{\Hh}=\mathfrak a(Sf,Sf)\ge0\), so \(S\)
is positive.  Finally, if \(Sf=0\), then \eqref{eq:state-weak} gives
\((f,v)_{\Hh}=0\) for all \(v\in\V\).  Since \(C_c^\infty(D)\) is dense in
\(L^2(D)\) and \(H\)-valued simple functions are dense in \(\Hh\), the space
\(\V\) is dense in \(\Hh\). Then \(f=0\), and \(S\) is injective.
\end{proof}
Substituting \(y=Su\) eliminates the state from \eqref{eq:ocp} and leaves the
reduced primal problem
\begin{equation}
\label{eq:primal}
 \min_{u\in\C}
 J(u)=\frac12\norm{Su-y_d}_{\Hh}^2
 +\frac\gamma2\norm{u}_{\Hh}^2,
\end{equation}
which is the form used throughout the paper.

\begin{proposition}[Existence and uniqueness]
\label{prop:primal-wp}
Problem \eqref{eq:primal} has a unique minimizer
\(u^*\){\(\in\C\)}.
Consequently, with the state equation in \eqref{eq:ocp} and the weak form \eqref{eq:state-weak}, \((Su^*,u^*)\in\V\times\C\) is the
unique solution of \eqref{eq:ocp}.
\end{proposition}
\begin{proof}
The extended functional
\[
 j(u)=\frac12\norm{Su-y_d}_{\Hh}^2+
 \frac\gamma2\norm{u}_{\Hh}^2+\I_\C(u),
\]
where \(\I_\C\) is the indicator function of \(\C\) (\(0\) on \(\C\)
and \(+\infty\) elsewhere),
and it is proper, coercive, convex, and weakly lower semicontinuous.  The direct
method gives a minimizer.  The term
\(\gamma\norm{u}_{\Hh}^2/2\) is strictly convex, so the minimizer is unique.
By Proposition~\ref{prop:state-operator}, the feasible pairs of
\eqref{eq:ocp} are exactly \((Su,u)\) with \(u\in\C\), and the objective of
\eqref{eq:ocp} at such a pair equals \(J(u)\); this gives the second
assertion.
\end{proof}

\section{Fenchel dual and the optimality system}\label{sec:fenchel}

The reduced problem \eqref{eq:primal} is a constrained minimization over
\(\C\).  This section removes the constraint by dualization.  We first compute
the two convex conjugates, then invoke Fenchel--Rockafellar duality to obtain a
dual problem in a residual variable, and finally change variables to the
elliptic adjoint.  The last form is the one the algorithm of
Section~\ref{sec:ssn} actually solves, and the recovery relations obtained at
the end of the section are what return the state and the control once the
adjoint is known.

\subsection{The dual problem}

To apply Fenchel--Rockafellar duality, we define
\(F:\Hh\to(-\infty,\infty]\) and \(G:\Hh\to\R\) by
\[
 F(u)=\frac\gamma2\norm{u}_{\Hh}^2+\I_\C(u),\qquad
 G(z)=\frac12\norm{z-y_d}_{\Hh}^2.
\]
Problem \eqref{eq:primal} can then be written as the unconstrained problem
\[
 \inf_{u\in\Hh}\{F(u)+G(Su)\}.
\]
The box constraint has been absorbed into \(F\) through the indicator
\(\I_\C\), so that both terms are now extended-real-valued convex functions on
\(\Hh\) and no explicit feasible set remains.  The next lemma shows that
conjugation turns the indicator into a finite, differentiable function: this is
the precise sense in which duality trades a constraint for smoothness.

\begin{lemma}[Conjugates]
\label{lem:conjugates}
For \(q,v\in\Hh\),
\[
 G^*(q)=\frac12\norm{q}_{\Hh}^2+(q,y_d)_{\Hh},
\qquad
 F^*(v)=
 \int_{\Gamma_\xi}\int_D\psi_\gamma(v(x,\xi))\,dx\,d\rho,
\]
where
\begin{equation}
\label{eq:psi}
 \psi_\gamma(s)=
 \begin{cases}
 as-\dfrac\gamma2a^2,&s<\gamma a,\\[1mm]
 \dfrac{s^2}{2\gamma},&\gamma a\le s\le\gamma b,\\[1mm]
 bs-\dfrac\gamma2b^2,&s>\gamma b.
 \end{cases}
\end{equation}
Moreover,
\begin{equation}
\label{eq:psi-derivative}
 \psi_\gamma'(s)=\Proj_{[a,b]}(s/\gamma).
\end{equation}
\end{lemma}

\begin{proof}
For \(G^*\), the substitution \(w=z-y_d\) gives
\[
 \sup_z\{(q,z)_{\Hh}-\tfrac12\norm{z-y_d}_{\Hh}^2\}
 =(q,y_d)_{\Hh}
 +\sup_w\{(q,w)_{\Hh}-\tfrac12\norm{w}_{\Hh}^2\}.
\]
For \(F^*\), maximize the scalar function
\[
 sr-\frac\gamma2r^2
\]
over \(r\in[a,b]\).  
Its maximizer is \(\Proj_{[a,b]}(s/\gamma)\), which gives \eqref{eq:psi}. Pointwise integration is justified by normal-integrand conjugacy; it also follows directly because the measurable pointwise maximizer belongs to \(\C\).
\end{proof}

Note that \(\psi_\gamma\) is continuously differentiable with a piecewise
linear derivative \eqref{eq:psi-derivative}: the nonsmoothness of the primal
constraint has been reduced to a single projection, which is exactly the
structure the semismooth Newton method of Section~\ref{sec:ssn} exploits.

Strong duality needs a qualification condition.  We use the following classical
form, in which the qualification is the finiteness and continuity of \(G\) at a single point
of the image of \(\dom F\); see~\cite{MR1727362,MR274683,Bauschke-2017-CAM}
for proofs and for more general statements.

\begin{theorem}[Fenchel--Rockafellar~{\cite[Chap.~III]{MR1727362}, \cite[Chap.~15]{Bauschke-2017-CAM}}]
\label{thm:fenchel-rockafellar}
Let \(X,Y\) be Banach spaces, let \(T:X\to Y\) be bounded and linear, and
let \(F:X\to(-\infty,\infty]\), \(G:Y\to(-\infty,\infty]\) be proper and
convex.  If there exists
\(x_0\in\dom F\) such that \(G\) is finite and continuous at \(Tx_0\), then
\[
 \inf_{x\in X}\{F(x)+G(Tx)\}
 =
 \mathop{\sup}\limits_{q\in Y'}\{-F^*(-T^*q)-G^*(q)\}.
\]
When the primal value is finite, the {supremum} is attained.
\end{theorem}
The tracking term \(G\) is a finite quadratic, hence continuous on all of
\(\Hh\), and \(\dom F=\C\) is nonempty because the constant control
\((a+b)/2\) is admissible.  No interior-point or Slater-type condition on the
box is needed. 

\begin{theorem}[Residual-variable Fenchel dual]
\label{thm:q-dual}
Under \eqref{eq:UE}, the Fenchel dual of \eqref{eq:primal} is
\begin{equation}
\label{eq:q-dual}
 \max_{q\in\Hh}
 \left\{
 -\int_{\Gamma_\xi\times D}\psi_\gamma(-Sq)\,d\rho\,dx
 -\frac12\norm{q}_{\Hh}^2
 -(q,y_d)_{\Hh}
 \right\}.
\end{equation}
There is no duality gap, and the maximum is attained at a unique
\(q^*\in\Hh\).
\end{theorem}

\begin{proof}
Apply Theorem~\ref{thm:fenchel-rockafellar} with
\[
 X=Y=\Hh,\qquad T=S,
\]
identifying \(\Hh'\) with \(\Hh\) through the Riesz map, so that
\(T^*=S^*\) is the Hilbert-space adjoint.
The quadratic \(G\) is finite and continuous everywhere, and the constant
midpoint control \((a+b)/2\) belongs to \(\dom F=\C\).  Thus the
qualification holds.  Since \(S^*=S\), Lemma~\ref{lem:conjugates} gives
\eqref{eq:q-dual}.
The map \(q\mapsto F^*(-Sq)\) is finite and convex, and
\(-\frac12\norm{q}_{\Hh}^2-(q,y_d)_{\Hh}\) is \(1\)-strongly concave. Hence the
dual objective is \(1\)-strongly concave and its maximizer is
unique.
\end{proof}

\subsection{The elliptic adjoint form and primal--dual recovery}

Problem \eqref{eq:q-dual} is already unconstrained, but it is not yet a
convenient computational form: the solution operator \(S\) sits inside the
nonlinearity \(\psi_\gamma(-Sq)\), so every evaluation of the objective
requires an elliptic solve.  Substituting \(p=-Sq\) moves \(S\) out of the
nonlinear term and replaces it by the differential operator \(A\) in the
quadratic term, which is what makes the discretized problem sparse and the
Jacobian applications of Section~\ref{sec:ssn} matrix-free.  For this
substitution to be a genuine reparametrization---creating no spurious
maximizers, losing none of the existing ones, and not tacitly shrinking the
feasible set---we must know that \(p\mapsto-Ap\) maps \(D(A)\) bijectively onto
\(\Hh\).  Establishing that is the sole purpose of the next lemma, and it is
the only place where the graph domain \eqref{eq:DA} is needed.

Let \(A\) denote the elliptic operator associated with \(\mathfrak a\).
As an unbounded operator on \(\Hh\), its natural graph domain is
\begin{equation}
\label{eq:DA}
 D(A)=\left\{p\in\V:Ap\in\Hh\right\}\,
\end{equation}
where \(Ap\in\Hh\) means that there is an element of \(\Hh\), denoted
\(Ap\), with
\[
 \mathfrak a(p,v)=(Ap,v)_{\Hh}\qquad(v\in\V).
\]
This element is unique since \(\V\) is dense in \(\Hh\) (see the proof
of Proposition~\ref{prop:state-operator}).

\begin{lemma}
\label{lem:graph-isomorphism}
The map \(A:D(A)\to\Hh\) is a bijection with inverse \(S\).
Consequently \(p\mapsto q=-Ap\) is a bijection from \(D(A)\) onto
\(\Hh\).
\end{lemma}

\begin{proof}
For every \(f\in\Hh\), Proposition~\ref{prop:state-operator} gives a unique
\(p=Sf\in\V\) satisfying \(Ap=f\).  Hence \(p\in D(A)\), so \(A\) is
surjective and \(S\) is its right inverse.  If \(Ap=0\), coercivity of
\(\mathfrak a\) gives \(p=0\), so \(A\) is injective.  No \(H^2(D)\)
regularity is used or asserted.
\end{proof}

\begin{theorem}[Elliptic-adjoint PDE dual]
\label{thm:p-dual}
The residual dual \eqref{eq:q-dual} is equivalent, up to the constant
\(-\norm{y_d}_{\Hh}^2/2\), to
\begin{equation}
\label{eq:p-dual}
 \min_{p\in D(A)}
 \left\{
 \frac12\norm{Ap-y_d}_{\Hh}^2+
 \int_{\Gamma_\xi\times D}\psi_\gamma(p)\,d\rho\,dx
 \right\}.
\end{equation}
\end{theorem}

\begin{proof}
Let 
\[
 p=-Sq.
\]
Lemma~\ref{lem:graph-isomorphism} gives \(q=-Ap\), while
\(-Sq=p\).  Substitution into \eqref{eq:q-dual} gives
\[
 \max_{p\in D(A)}
 \left\{
 -\int\psi_\gamma(p)
 -\frac12\norm{Ap}_{\Hh}^2
 +(Ap,y_d)_{\Hh}
 \right\}.
\]
Changing the sign of the objective and completing the square yields
\eqref{eq:p-dual}.
\end{proof}

Both duals determine the adjoint, but the quantities of interest---and the
residuals used to monitor convergence---are the primal pair \((y,u)\).  The
final result of this section therefore records how the primal variables are
recovered from a dual solution and verifies that the recovered triple solves
the classical optimality system, so that nothing is lost in passing to the
dual.  The projection formula it produces is not merely a postprocessing step:
it is the term that reappears in the discrete equation of
Section~\ref{sec:ssn} and that the semismooth Newton method differentiates.

\begin{theorem}[Optimality system]
\label{thm:kkt}
Let \(u^*,y^*,q^*,p^*\) be primal and dual optimizers, with
\(y^*=Su^*\) and \(p^*=-Sq^*\).  Then
\[
 q^*=y^*-y_d,\qquad
 u^*=\Proj_{[a,b]}(p^*/\gamma),
\]
and the primal and dual problems are equivalent to
\begin{equation}
\label{eq:kkt}
\begin{cases}
 -\nabla\cdot(\kappa\nabla y)=u,
       &y|_{\partial D}=0,\\
 -\nabla\cdot(\kappa\nabla p)=y_d-y,
       &p|_{\partial D}=0,\\
 u=\Proj_{[a,b]}(p/\gamma).
\end{cases}
\end{equation}
Here \(p\) is the negative of the adjoint state in the usual Lagrangian
convention and the recovery formula followed as
\(u=\Proj_{[a,b]}(p/\gamma)\).
\end{theorem}

\begin{proof}
Fenchel equality is equivalent to
\[
 q^*\in\partial G(Su^*),\qquad
 -S^*q^*\in\partial F(u^*).
\]
The first inclusion gives
\[
 q^*=Su^*-y_d=y^*-y_d.
\]
Since \(S^*=S\) and \(p^*=-Sq^*\), the second inclusion is
\[
 p^*\in\gamma u^*+N_\C(u^*),
\]
here
\(N_\C(u)=\{v\in\Hh:(v,w-u)_{\Hh}\le0\text{ for all }w\in\C\}\) is the
normal cone of \(\C\) at \(u\in\C\), 
which is equivalent to
\[
 u^*=\Proj_\C(p^*/\gamma)
 =\Proj_{[a,b]}(p^*/\gamma)
\quad\text{a.e.}
\]
Finally,
\[
 Ap^*=-q^*=y_d-y^*.
\]
Together with primal feasibility, these identities give \eqref{eq:kkt}.
\end{proof}

\section{Stochastic Galerkin discretization and its dual}\label{sec:sgfe}

\subsection{SG--FE state and mass matrices}
Let \(\mathcal P_P\subset L^2_\rho(\Gamma_\xi)\) be the space of polynomials in \(\xi\) of total degree at most \(P\), and let \(\{\psi_i\}_{i=1}^{n_P}\), \(n_P=\binom{K+P}{P}\), be a basis of \(\mathcal P_P\) that is orthonormal in \(L^2_\rho(\Gamma_\xi)\), i.e., \(\E[\psi_i\psi_j]=\delta_{ij}\), which is the generalized polynomial chaos (gPC)  basis~\cite{xiu2003modeling}. Let \(\{\phi_r\}_{r=1}^{n_h}\) be the nodal basis of a conforming finite-element space \(V_h\subset H_0^1(D)\).
The SG method approximates the state in the tensor-product space $\mathcal P_P\otimes V_h$,
\[
 y_{hP}(x,\xi)=\sum_{i=1}^{n_P}\sum_{r=1}^{n_h}y_{(i,r)}\,\psi_i(\xi)\phi_r(x),
\]
and requires the weak state equation \eqref{eq:state-weak} to hold for all test functions $v=\psi_j\phi_s$, $j=1,\dots,n_P$, $s=1,\dots,n_h$. The coefficient vector $\mathbf y=(y_{(i,r)})\in\R^N$, $N\coloneqq n_Pn_h$, is ordered stochastic-major, i.e., the index $(i,r)$ is stored at position $(i-1)n_h+r$. The adjoint $\mathbf p$ is expanded in the same basis. 

Unlike the state, the control is not expanded in the gPC basis but represented by its values $\mathbf u=(u_{(k,r)})\in\R^m$, $m\coloneqq n_qn_h$, at the points $(x_r,\xi_k)$. This keeps the box constraint $a\le u\le b$ pointwise and componentwise. The SG load $\E[\int_D u\,\psi_j\phi_s\,dx]$ is approximated by the quadrature rule in $\xi$ and the lumped mass in $x$, which gives $(E^\top W\mathbf u)_{(j,s)}=\sum_{k}w_k\psi_j(\xi_k)(M_L)_{ss}u_{(k,s)}$.

Suppose that the diffusion coefficient $\kappa$ used in the SG form admits the separable representation
\[
 \kappa(x,\xi)=
 \sum_{\ell=0}^{N_K-1}\kappa_\ell(x)\zeta_\ell(\xi),
 \qquad \zeta_0\equiv1,
\]
where $\kappa_\ell\in L^\infty(D)$ are deterministic spatial functions, $\zeta_\ell\colon\Gamma_\xi\to\R$ are given scalar functions of the parameter, and $N_K$ is the number of terms. For the affine coefficients used in Section~\ref{sec:numerics}, $N_K=K+1$ and $\zeta_\ell(\xi)=\xi_\ell$ for $\ell\ge1$.
Followed by~\cite{Benner-2016-BDP}, with the definition of
\begin{align*}
 (G_\ell)_{ij}
 &=\E[\zeta_\ell\psi_i\psi_j],\\ (M_D)_{ij}
 &=\int_D\phi_i\phi_j\,dx,\\
 (K_\ell)_{ij}
 &=\int_D\kappa_\ell(x)\nabla\phi_i\cdot\nabla\phi_j\,dx,
\end{align*}
the SG mass matrix and stiffness matrix are
\begin{equation}
\label{eq:MK}
 \mathcal{M}=G_0\otimes M_D,\qquad
 \mathcal{K}=\sum_{\ell=0}^{N_K-1}G_\ell\otimes K_\ell.
\end{equation}
For an exactly integrated orthonormal basis, \(G_0=I\).  
The control is represented by its values at the points
$(x_r,\xi_k)$, where $x_1,\ldots,x_{n_h}$ are the interior nodes of $V_h$
and $\{{\xi_k,w_k}\}_{k=1}^{n_q}$ is a positive quadrature rule on
$\Gamma_\xi$.   Let
$V_\xi\in\R^{n_q\times n_P}$ with $(V_\xi)_{kj}={\psi_j(\xi_k)}$,
$W_\xi=\operatorname{diag}(w_1,\ldots,w_{n_q})$, and
$M_L=\operatorname{diag}(M_D\mathbf 1)$, the row-sum lumped spatial mass.  
Define
\[
 E=V_\xi\otimes I_{n_h}\in\R^{n_qn_h\times n_Pn_h},\qquad
 W=W_\xi\otimes M_L\in\R^{n_qn_h\times n_qn_h}.
\]
For a coefficient vector $\mathbf p$, $(E\mathbf p)_{(k,r)}=\sum_{j}\psi_j(\xi_k)p_{(j,r)}$
is the value at $(x_r,\xi_k)$ of the SG function
$\sum_{j,r}p_{(j,r)}\psi_j(\xi)\phi_r(x)$, because the nodal basis satisfies
$\phi_s(x_r)=\delta_{rs}$.    The entries of $W$ are
positive, so $W$ is positive definite and diagonal.  Boundary nodes carry no control values,
since the Dirichlet degrees of freedom are eliminated.
Thus $E$ and $W$
belong to the control discretization, whereas $M_D$, $K_\ell$, $\mathcal M$,
and $\mathcal K$ belong to the state discretization.

We assume that \(\mathcal M\) and \(\mathcal K\) are symmetric positive
definite; this holds, for instance, when the matrices \(G_\ell\) are computed
exactly and $\kappa$ satisfies \eqref{eq:UE}.


\subsection{Lagrange dual and KKT system}
Let  \(\mathbf u\in\R^m\) be independent control values at the positive quadrature
points. Observe that the discrete state $\mathbf y\in\R^N$ satisfies the SG weak equation
\[
 \Kcal \mathbf y=E^\top W\mathbf u.
\]
The exact finite-dimensional constrained primal is
\begin{equation}
\label{eq:discrete-primal}
\begin{aligned}
 \min_{\mathbf y\in\R^N,\ \mathbf u\in\R^m}\quad&
 \frac12(\mathbf y-\mathbf y_d)^\top \Mcal (\mathbf y-\mathbf y_d)
 +\frac\gamma2\mathbf u^\top W\mathbf u,\\
 \text{subject to}\quad&
 \Kcal \mathbf y-E^\top W\mathbf u=0,\\
 &a\le \mathbf u\le b.
\end{aligned}
\end{equation}
\cref{thm:discrete-dual} below refers to the dual of \eqref{eq:discrete-primal}.  

\begin{theorem}
\label{thm:discrete-dual}
Problem
\eqref{eq:discrete-primal} has the concave dual
\begin{equation}
\label{eq:discrete-dual}
 \max_{\mathbf p\in\R^N}
 \left\{
 -\frac12\mathbf p^\top \Kcal \Mcal ^{-1}\Kcal^\top \mathbf p
 +\mathbf p^\top \Kcal \mathbf y_d
 -\sum_{\ell=1}^mw_\ell\psi_\gamma((E\mathbf p)_\ell)
 \right\}.
\end{equation}
There is no duality gap.  Moreover, the primal and dual optimizers are unique.
\end{theorem}

\begin{proof}
Let
\[
 \bar{\mathbf u}=\frac{a+b}{2}\mathbf 1,\qquad
 \bar{\mathbf y}=\Kcal^{-1}E^\top W\bar{\mathbf u}.
\]
This is an affine-feasible point in the relative interior of the box, so
the finite-dimensional relative-interior qualification holds.

Now, consider the Lagrangian
\begin{equation}
\label{lagrange}
 \mathcal L_h(\mathbf y,\mathbf u,\mathbf p)=
 \frac12(\mathbf y-\mathbf y_d)^\top \Mcal(\mathbf y-\mathbf y_d)
 +\frac\gamma2\mathbf u^\top W\mathbf u
 +\mathbf p^\top(\Kcal \mathbf y-E^\top W\mathbf u)
 +\I_{[a,b]^m}(\mathbf u).
 \end{equation}
Minimization in \(\mathbf y\) of \eqref{lagrange} gives
\[
 \mathbf y=\mathbf y_d-\Mcal^{-1}\Kcal ^\top \mathbf p
\]
and the minimum objective
\[
 -\frac12\mathbf p^\top \Kcal \Mcal^{-1}\Kcal^\top \mathbf p+\mathbf p^\top \Kcal \mathbf y_d.
\]
For \(s_\ell=(E\mathbf p)_\ell\),
\[
 \inf_{r\in[a,b]}
 w_\ell\left(\frac\gamma2r^2-s_\ell r\right)
 =-w_\ell\psi_\gamma(s_\ell).
\]
Since $W=\operatorname{diag}(w_1,\dots,w_m)$ is diagonal and the box is a product of intervals, the $\mathbf u$-part of $\mathcal L_h$ separates into the $m$ scalar problems above. Summing them over $\ell=1,\dots,m$ and adding the minimum in $\mathbf y$ gives \eqref{eq:discrete-dual}.
Since \(\Kcal\) is invertible,
\(\Kcal \Mcal^{-1} \Kcal^\top\succ0\), so the dual objective is {strongly} concave and
tends to \(-\infty\) as \(\norm{\mathbf p}\to\infty\).  Hence it has a unique
maximizer.  The primal objective is strictly convex on the affine feasible
set, giving primal uniqueness.
\end{proof}

\begin{corollary}
\label{cor:two-field}
After eliminating \(\mathbf u\), the KKT conditions are
\begin{subequations}
\label{eq:two-field}
\begin{align}
 \Mcal \mathbf y+\Kcal^\top \mathbf p&=\Mcal \mathbf y_d, \label{eq:two-field-a}\\
 \Kcal \mathbf y&=R(\mathbf p), \label{eq:two-field-b}
\end{align}
\end{subequations}
where \(R(\mathbf p)=E^\top W\Proj_{[a,b]}(E\mathbf p/\gamma)\).
\end{corollary}
\begin{proof}
Stationarity in \(\mathbf y\), \(\nabla_{\mathbf y}\mathcal L_h=\Mcal(\mathbf y-\mathbf y_d)+\Kcal^\top \mathbf p=0\), gives \eqref{eq:two-field-a}. Stationarity in each
control value is
\[
 0\in\gamma \mathbf u-E\mathbf p+N_{[a,b]^m}(\mathbf u),
\]
where \(N_{[a,b]^m}\) is the normal cone of the box \([a,b]^m\subset\R^m\),
because positive diagonal rescaling leaves the componentwise box normal
cone unchanged.  Therefore
\[
 \mathbf u=\Proj_{[a,b]}(E\mathbf p/\gamma).
\]
Substitution into the primal feasibility condition \(\Kcal \mathbf y-E^\top W\mathbf u=0\) gives \eqref{eq:two-field-b}.  Equivalently,
differentiating \eqref{eq:discrete-dual} and using
\eqref{eq:psi-derivative} gives
\[
 -\Kcal \Mcal^{-1}\Kcal^\top \mathbf p+\Kcal \mathbf y_d-R(\mathbf p)=0,
\]
which is the same equation after substituting
\(\mathbf y=\mathbf y_d-\Mcal^{-1}\Kcal^\top \mathbf p\).
\end{proof}


\begin{lemma}
\label{lem:monotone}
If \(W\succ0\), then \(R:\R^N\to\R^N\) in \eqref{eq:two-field} is
monotone and Lipschitz continuous.
\end{lemma}

\begin{proof}
For any \(\mathbf v,\mathbf z\in\R^N\),
\begin{align*}
 &(\mathbf v-\mathbf z)^\top(R(\mathbf v)-R(\mathbf z))=(E \mathbf v-E \mathbf z)^\top W
 \left[
 \Proj_{[a,b]}(E\mathbf v/\gamma)
 -\Proj_{[a,b]}(E\mathbf z/\gamma)
 \right]\ge0.
\end{align*}
Each summand is nonnegative because the scalar projection is monotone and
each quadrature weight is positive.  The scalar projection is
\(1\)-Lipschitz, so \(R\) is Lipschitz.  
\end{proof}

\section{A semismooth Newton--PCG method}
\label{sec:ssn}

Eliminating $\mathbf y$ in \eqref{eq:two-field} gives
\begin{equation}\label{eq:F}
F(\mathbf p)=\Kcal \Mcal ^{-1}\Kcal ^\top \mathbf p-\Kcal \mathbf y_d+E^\top W\Proj_{[a,b]}(E\mathbf p/\gamma)=0.
\end{equation}
For the diagonal free-set mask
$D(\mathbf p)=\operatorname{diag}(d_q(\mathbf p))\in\mathbb R^{m\times m}$, where
$d_q(\mathbf p)=1$ on a free component and $d_q(\mathbf p)=0$ on a saturated component
(and zero is selected at a breakpoint), the Newton derivative of $F$ in the sense
of~\cite{Hintermuller-2002-PDA,Ulbrich-2011-SNM} is
\begin{equation}\label{eq:J}
J_D=\Kcal \Mcal^{-1}\Kcal^\top+\gamma^{-1}E^\top W D(\mathbf p)E.
\end{equation}
For a vector $\mathbf v$, $J_D\mathbf v$ is evaluated matrix-free using one application each of $\Kcal$, $\Kcal^\top$, $\Mcal^{-1}$, $E$, and $E^\top$. Each Newton system is solved by PCG to the {forcing} tolerance
\begin{equation}\label{eq:ew}
\eta_k=\min\{{\eta_{\max}},\max(10^{-12},0.1\sqrt{r_k})\},\qquad r_k=\norm{F(\mathbf p_k)}{/\max\{1,\norm{\Kcal \mathbf y_d}\}},
\end{equation}
and the step is globalized by Armijo backtracking with parameters $10^{-4}$
and $1/2$, which preserves descent of the convex dual merit function.  
The
iteration stops when the normalized KKT residual is at most {$10^{-6}$} or after
{$k_{\max}$} Newton steps, and the values of $\eta_{\max}$ and $k_{\max}$ are given in
Section~\ref{sec:numerics}.  The complete procedure is specified by
Algorithm~\ref{alg:ssn-pcg} in Section~\ref{sec:numerics}, after the
preconditioner has been specified.

\subsection{Current-set preconditioner and spectral bounds}
\label{sec:woodbury}
Here $M_D$, $K_0$, and $K_\ell$ retain the spatial finite-element
meanings defined in Section~\ref{sec:sgfe}.  The calligraphic symbols
$\mathcal{M}$, $\mathcal{K}$, and $\overline{\mathcal{K}}$ denote the full
stochastic--Galerkin matrices.
We use stochastic-major ordering first and apply an orthogonal
permutation to node-major ordering for the local assembly estimates.  
In stochastic-major
ordering, since $G_0=I$ for the exactly integrated orthonormal basis of Section~\ref{sec:sgfe}, the matrices in \eqref{eq:MK}, with the mean term $\ell=0$ written separately, become
\begin{equation}\label{eq:wm-SG}
 \mathcal{M}=I_{n_P}\otimes M_D,\qquad
 \mathcal{K}=I_{n_P}\otimes K_0+
       \sum_{\ell=1}^{N_K-1}G_\ell\otimes K_\ell,
 \qquad \overline{\mathcal{K}}=I_{n_P}\otimes K_0.
\end{equation}
Here $\overline{\mathcal K}$ is the part of $\mathcal K$ generated by the mean coefficient $\kappa_0$, $N_K-1$ is the number of non-mean terms of $\kappa$, and $M_D,K_0$ are symmetric positive definite, $G_\ell,K_\ell$
are symmetric, and $\mathcal{K}$ is positive definite. The stochastic basis is
orthonormal for the mass inner product represented in $\mathcal{M}$. All coefficient
amplitudes and normalization factors are already included in $G_\ell$
or $K_\ell$. 

The control weights $W$ are positive and diagonal, as above.  Written in
terms of the reduced equation \eqref{eq:F}, the Jacobian matrix \eqref{eq:J}
associated with a selected projection derivative $D$ has the form
\begin{equation}\label{eq:wm-JP}
 H=\mathcal{K}\mathcal{M}^{-1}\mathcal{K}^\top,\quad H_0=\overline{\mathcal{K}}\mathcal{M}^{-1}\overline{\mathcal{K}}^\top,
 \qquad J=H+\gamma^{-1}C_D,\quad P_W=H_0+\gamma^{-1}C_D,
\end{equation}
where $\gamma>0$ and $C_D=E^\top WDE\succeq0$. 
The algorithm uses the binary mask $D=D(\mathbf p)$ of \eqref{eq:J} and the spectral results below hold for any positive semidefinite $C_D$ shared by $J$ and $P_W$.

Let $\|\mathbf v\|_Q=(\mathbf v^\top Q\mathbf v)^{1/2}$ for an SPD matrix $Q$.
All unweighted matrix norms below are Euclidean operator norms.
The notation $\rm cond_ {\rm sp}(P_W^{-1}J)$ means the ratio of extreme
eigenvalues, equivalently the condition number of the SPD matrix
$P_W^{-1/2}JP_W^{-1/2}$.

\begin{lemma}
\label{lem:wm-reduction}
Let
\[
 \mathcal T=\mathcal{M}^{-1/2}\mathcal{K}\overline{\mathcal{K}}^{-1}\mathcal{M}^{1/2},\qquad
 \sigma_-=\sigma_{\min}(\mathcal T)>0,\quad \sigma_+=\sigma_{\max}(\mathcal T),
\]
where $\sigma_{\min}(\mathcal T)$ and $\sigma_{\max}(\mathcal T)$ denote the smallest and largest singular values of~$\mathcal T$, 
and the extreme generalized eigenvalues of $(H,H_0)$ are exactly
$\sigma_-^2,\sigma_+^2$. In particular,
\begin{equation}\label{eq:wm-Hbounds}
 \sigma_-^2H_0\preceq H\preceq \sigma_+^2H_0.
\end{equation}
\end{lemma}
\begin{proof}
Let $L_0=\mathcal{M}^{-1/2}\overline{\mathcal{K}}$, an invertible matrix, since $\mathcal M\succ0$ and $\overline{\mathcal K}=I_{n_P}\otimes K_0$ is invertible. With the fact that
$\mathcal K=\mathcal K^\top$ and $\overline{\mathcal K}=\overline{\mathcal K}^\top$
from~\eqref{eq:wm-SG}, direct multiplication
gives $H_0=L_0^\top L_0$ and $H=L_0^\top \mathcal T^\top \mathcal TL_0$, indeed $\mathcal TL_0=\mathcal M^{-1/2}\mathcal K$. Therefore
\[
 \frac{\mathbf v^\top H\mathbf v}{\mathbf v^\top H_0\mathbf v}
 =\frac{\|\mathcal T\mathbf z\|_2^2}{\|\mathbf z\|_2^2},\qquad \mathbf z=L_0\mathbf v.
\]
Since $L_0$ is onto, extrema over $\mathbf v\ne0$ equal extrema over $\mathbf z\ne0$.
They are the extreme eigenvalues of $\mathcal T^\top \mathcal T$, namely $\sigma_-^2$ and $\sigma_+^2$. Note that $\sigma_->0$ because $\mathcal T$ is invertible.
\end{proof}
\begin{theorem}
\label{thm:wm-algebra}
If $cH_0\preceq H\preceq CH_0$ for constants $0<c\le C$, then
for every $\gamma>0$ and $C_D\succeq0$,
\begin{equation}\label{eq:wm-spectrum}
 \operatorname{spec}(P_W^{-1}J)
 \subset[\theta_-,\theta_+],\qquad
 \theta_-=\min\{c,1\},\quad\theta_+=\max\{C,1\}.
\end{equation}
Consequently $\rm cond_ {\rm sp}(P_W^{-1}J)\le\theta_+/\theta_-$.
In particular, one can always take $c=\sigma_-^2,C=\sigma_+^2$ at a fixed discretization.
\end{theorem}
\begin{proof}
First, observe that both
\[
 J-\theta_-P_W=(H-\theta_-H_0)
                   +(1-\theta_-)\gamma^{-1}C_D
\]
and
\[
 \theta_+P_W-J=(\theta_+H_0-H)
                   +(\theta_+-1)\gamma^{-1}C_D
\]
are positive semidefinite: in the first terms  because $\theta_-\le c$ and $\theta_+\ge C$, and in the second terms because $\theta_-\le1\le\theta_+$ and $C_D\succeq0$. 
Since $P_W\succ0$, congruence with $P_W^{-1/2}$ preserves semidefiniteness and gives $\theta_-I\preceq P_W^{-1/2}JP_W^{-1/2}\preceq\theta_+I$. Hence the eigenvalues of the symmetric matrix $P_W^{-1/2}JP_W^{-1/2}$ lie in $[\theta_-,\theta_+]$, and so do those of $P_W^{-1}J=P_W^{-1/2}(P_W^{-1/2}JP_W^{-1/2})P_W^{1/2}$.
\end{proof}


\subsection{Uniform bounds from local assembly}
We next provide sufficient local structure that yields a uniform
bound. 
Consider a square of fixed side length $L_D>0$ with $n_e\ge2$ equal
elements per direction, $h=L_D/n_e$, and conforming bilinear elements.
All boundary nodal values are zero and are eliminated. On each element
$T$, let $q$ index the four tensor Gauss points $x_{Tq}$, let $w_{Tq}>0$ include
the physical Jacobian, let $\mathbf g_{Tiq}\in\mathbb R^2$ be the gradient
of the local nodal basis function $i$ at $x_{Tq}$, and let $a_{\ell,Tq}=\kappa_\ell(x_{Tq})$.
\begin{assumption}\label{ass:wm-local}
Let the scalar reference weights and symmetric coefficient blocks satisfy
\begin{equation}\label{eq:wm-localbounds}
 0<\nu_-\le a_{0,Tq}\le\nu_+,
 \qquad \alpha I_{n_P}\preceq A_{Tq}\preceq\beta I_{n_P},
 \qquad 0<\alpha\le\beta.
\end{equation}
The reference stiffness has entries
\[
 (K_0)_{ij}=\sum_{T\ni i,j}\sum_q
                 w_{Tq}a_{0,Tq}\,\mathbf g_{Tiq}^\top \mathbf g_{Tjq}.
\]
After an orthogonal permutation into node-major ordering, the full
stiffness has $n_P\times n_P$ blocks
\begin{equation}\label{eq:wm-localK}
K_{ij}=\sum_{T\ni i,j}\sum_q
      w_{Tq}a_{0,Tq}(\mathbf g_{Tiq}^\top \mathbf g_{Tjq})A_{Tq}.
\end{equation}
The spatial mass is the consistent bilinear mass, integrated exactly
by the same Gauss rule. 
\end{assumption}

Orthogonal reordering does not change any spectral or norm assertion.
For a Kronecker assembly with scalar coefficient arrays $a_{\ell,Tq}$,
\eqref{eq:wm-localK} is exactly \eqref{eq:wm-SG} when
\begin{equation}\label{eq:wm-AfromG}
 A_{Tq}=I_{n_P}+
       \sum_{\ell=1}^{N_K-1}\frac{a_{\ell,Tq}}{a_{0,Tq}}G_\ell.
\end{equation}

\begin{lemma}
\label{lem:wm-reference}
Under Assumption~\ref{ass:wm-local},
\begin{equation}\label{eq:wm-referencebounds}
 \frac{h^2}{9}I\preceq M_D\preceq h^2I,\qquad
 (K_0)_{ii}\le\frac{8\nu_+}{3},\qquad
 K_0\succeq\lambda_*M_D,
 \quad\lambda_*\coloneqq\frac{8\nu_-}{L_D^2}.
\end{equation}
Moreover $\alpha\overline{\mathcal{K}}\preceq \mathcal{K}\preceq\beta\overline{\mathcal{K}}$.
\end{lemma}
\begin{proof}
Let one-dimensional interior matrices be
\[
 M_1=\frac h6\operatorname{tridiag}(1,4,1),\qquad
 S_1=\frac1h\operatorname{tridiag}(-1,2,-1).
\]
Inserting the vector with entries $\sin(ir\pi/n_e)$ into the three-point stencils shows that the eigenvalues are
\[
 m_r=\frac h3(2+\cos\theta_r),\qquad
 s_r=\frac2h(1-\cos\theta_r),\qquad\theta_r=\frac{r\pi}{n_e}.
\]
These orthogonal sine vectors form a basis. Since
$h/3\le m_r\le h$, the tensor identity $M_D=M_1\otimes M_1$
gives the mass bounds. The unit-reference stiffness is
$K_{\rm unit}=M_1\otimes S_1+S_1\otimes M_1$.
Its generalized eigenvalues relative to $M_D$ are
$s_r/m_r+s_t/m_t$. The ratio is increasing with $\theta_r$ and
\[
 \frac{s_1}{m_1}
 =\frac{12}{h^2}\frac{\sin^2(\pi/(2n_e))}{2+\cos(\pi/n_e)}
 \ge\frac{4}{h^2n_e^2}=\frac4{L_D^2}.
\]
Here $\sin x\ge2x/\pi$ for $0\le x\le\pi/2$ follows from
concavity of sine. Positive quadrature gives
$K_0\succeq\nu_-K_{\rm unit}$, proving the generalized lower bound.

The diagonal of the unit local $4\times4$ bilinear stiffness is
$2/3$. Indeed, for $N=(1-r)(1-s)/4$ on $[-1,1]^2$,
the physical gradient and Jacobian give the integral
$\frac1{16}\int_{-1}^1\int_{-1}^1((1-s)^2+(1-r)^2)\,dr\,ds=2/3$.
Gauss quadrature is exact here. An interior node belongs to four
elements, yielding $(K_0)_{ii}\le4\nu_+(2/3)$.

Finally let $\mathbf v_i\in\mathbb R^{n_P}$ be nodal values, extended by
zero at boundary nodes, and let
\[
 U_{Tq}(\mathbf v)=\sum_{i\in T}\mathbf v_i \mathbf g_{Tiq}^\top
                   \in\mathbb R^{n_P\times2},
\]
then
\[
 \mathbf v^\top\overline{\mathcal{K}} \mathbf v=\sum_{T,q}w_{Tq}a_{0,Tq}\|U_{Tq}(\mathbf v)\|_F^2,
\]
 while $\mathbf v^\top\mathcal{K}\mathbf v$ replaces the squared Frobenius norm by
$\operatorname{tr}(U_{Tq}(\mathbf v)^\top A_{Tq}U_{Tq}(\mathbf v))$.
The eigenvalue bounds on $A_{Tq}$ prove the final assertion.
\end{proof}

For each interior node $i$, let $\mathcal Q_i$ consist of all incident
element--quadrature pairs. We take $A_i$ to be the arithmetic average
\begin{equation}\label{eq:wm-average}
 A_i=\frac1{|\mathcal Q_i|}\sum_{(T,q)\in\mathcal Q_i}A_{Tq},
 \qquad
 \Omega_h=\max_i\max_{(T,q)\in\mathcal Q_i}
                   \frac{\|A_{Tq}-A_i\|_2}{h}.
\end{equation}
Each $A_i$ is a convex combination of the symmetric blocks $A_{Tq}$, and
the set $\{X=X^\top:\alpha I\preceq X\preceq\beta I\}$ is convex, then
\eqref{eq:wm-localbounds} gives $\alpha I\preceq A_i\preceq\beta I$.
Consequently the matrix $\mathsf D$ defined next is symmetric positive
definite, with $\|\mathsf D\|_2\le\beta$ and
$\|\mathsf D^{-1}\|_2\le\alpha^{-1}$.
Define the nodal block-diagonal multiplication matrix
$\mathsf D=\operatorname{diag}_i(A_i)$ and set
\begin{equation}\label{eq:wm-decomposition}
 \mathcal{K}=\mathsf D\overline{\mathcal{K}}+R.
\end{equation}

\begin{lemma}\label{lem:wm-remainder}
Under Assumption~\ref{ass:wm-local}, the remainder in \eqref{eq:wm-decomposition} satisfies
\begin{equation}\label{eq:wm-Rbound}
 \|R\mathbf v\|_{\mathcal{M}^{-1}}\le C_R\Omega_h\|\mathbf v\|_{\overline{\mathcal{K}}},
 \qquad C_R=\sqrt{96\nu_+}.
\end{equation}
\end{lemma}
\begin{proof}
Subtracting the two assembled rows before the boundary values are eliminated, we find that the $i$th block row is exactly
\[
 (R\mathbf v)_i=\sum_{(T,q)\in\mathcal Q_i}
 w_{Tq}a_{0,Tq}(A_{Tq}-A_i)U_{Tq}(\mathbf v)\mathbf g_{Tiq}.
\]
Since $\|U g\|_2\le\|U\|_F\|g\|_2$, weighted Cauchy--Schwarz gives
\begin{align*}
 \|(R\mathbf v)_i\|_2^2
 &\le h^2\Omega_h^2
 \left(\sum_{\mathcal Q_i}w_{Tq}a_{0,Tq}\|\mathbf g_{Tiq}\|_2^2\right)
 \left(\sum_{\mathcal Q_i}w_{Tq}a_{0,Tq}\|U_{Tq}(\mathbf v)\|_F^2\right).
\end{align*}
The first parenthesis equals $(K_0)_{ii}\le8\nu_+/3$.
Upon summing over interior nodes, any element energy is counted at
most four times. 
Thus
\[
 \|R\mathbf v\|_2^2\le\frac{32\nu_+}{3}h^2\Omega_h^2
                            \mathbf v^\top\overline{\mathcal{K}}\mathbf v.
\]
By \eqref{eq:wm-referencebounds}, $M_D^{-1}\preceq9h^{-2}I$, hence $\mathcal M^{-1}=I_{n_P}\otimes M_D^{-1}\preceq9h^{-2}I$ and
$$\|R \mathbf v\|_{\mathcal M^{-1}}^2=(R \mathbf v)^\top\mathcal M^{-1}(R \mathbf v)\le(R \mathbf v)^\top(9h^{-2}I)(R \mathbf v)=9h^{-2}\|R \mathbf v\|_2^2  \le96\,\nu_+\Omega_h^2\,\|\mathbf v\|_{\overline{\mathcal K}}^2 .
$$
Taking square roots gives \eqref{eq:wm-Rbound} with $C_R=\sqrt{96\nu_+}$.
All vector bounds use operator norms of the coefficient blocks.
The Frobenius norm above sums the physical gradient energies, not the
entries of $A_{Tq}$.
\end{proof}
\begin{lemma}[Multiplication and discrete energy estimates]
\label{lem:wm-multiplication}
For all $\mathbf f\in\R^N$,
\begin{equation}\label{eq:wm-Dnorms}
 \|\mathsf D\mathbf f\|_{\mathcal{M}^{-1}}\le3\beta\|\mathbf f\|_{\mathcal{M}^{-1}},\qquad
 \|\mathsf D^{-1}\mathbf f\|_{\mathcal{M}^{-1}}\le\frac3\alpha\|\mathbf f\|_{\mathcal{M}^{-1}}.
\end{equation}
With $c_*\coloneqq\lambda_*^{-1/2}$, one also has
\begin{equation}\label{eq:wm-energy}
 \|\mathbf v\|_{\overline{\mathcal{K}}}\le c_*\|\overline{\mathcal{K}}\mathbf v\|_{\mathcal{M}^{-1}},\qquad
 \|\mathbf v\|_{\overline{\mathcal{K}}}\le\frac{c_*}{\alpha}\|\mathcal{K}\mathbf v\|_{\mathcal{M}^{-1}}.
\end{equation}
\end{lemma}
\begin{proof}
Since the spectral condition number $\operatorname{cond}_2(\mathcal M)=\lambda_{\max}(\mathcal M)/\lambda_{\min}(\mathcal M)$,
the operator norm of the nodal multiplication, 
which multiplies the gPC coefficients at each node $i$ by $A_i$, satisfies
\[
 \|\mathcal{M}^{-1/2}\mathsf D\mathcal{M}^{1/2}\|_2
 \le\sqrt{\operatorname{cond}_2(\mathcal{M})}\|\mathsf D\|_2\le3\beta.
\]
Applying the same argument to $\mathsf D^{-1}$ yields $3/\alpha$.

Since $\overline{\mathcal K}=I_{n_P}\otimes K_0\succeq\lambda_*\mathcal M$ by \eqref{eq:wm-referencebounds}, we have $\|\mathbf v\|_{\mathcal M}\le c_*\|\mathbf v\|_{\overline{\mathcal K}}$. The weighted Cauchy--Schwarz inequality $|\mathbf a^\top\mathbf b|\le\|\mathbf a\|_{\mathcal M}\|\mathbf b\|_{\mathcal M^{-1}}$ with $\mathbf a=\mathbf v$, $\mathbf b=\overline{\mathcal K}\mathbf v$ then gives
\[
 \|\mathbf v\|_{\overline{\mathcal K}}^2=\mathbf v^\top\overline{\mathcal K}\mathbf v
 \le\|\mathbf v\|_{\mathcal M}\|\overline{\mathcal K}\mathbf v\|_{\mathcal M^{-1}}
 \le c_*\|\mathbf v\|_{\overline{\mathcal K}}\|\overline{\mathcal K}\mathbf v\|_{\mathcal M^{-1}},
\]
Similarly,
\[
 \alpha\|\mathbf v\|_{\overline{\mathcal{K}}}^2\le \mathbf v^\top \mathcal{K}\mathbf v
 \le\|\mathbf v\|_{\mathcal{M}}\|\mathcal{K}\mathbf v\|_{\mathcal{M}^{-1}}
 \le c_*\|\mathbf v\|_{\overline{\mathcal{K}}}\|\mathcal{K}\mathbf v\|_{\mathcal{M}^{-1}}.
\]
Dividing by $\|\mathbf v\|_{\overline{\mathcal K}}$ (the case $\mathbf v=0$ is trivial) yields \eqref{eq:wm-energy}.
\end{proof}


\begin{theorem}
\label{thm:wm-graph}
Under Assumption~\ref{ass:wm-local}, let
\begin{equation}\label{eq:wm-UV}
 q_h=C_Rc_*\Omega_h,\qquad
 U_h=3\beta+q_h,\qquad
 V_h=\frac3\alpha\left(1+\frac{q_h}{\alpha}\right).
\end{equation}
Then
\begin{equation}\label{eq:wm-graph}
 V_h^{-1}\|\overline{\mathcal{K}}\mathbf v\|_{\mathcal{M}^{-1}}
 \le\|\mathcal{K}\mathbf v\|_{\mathcal{M}^{-1}}
 \le U_h\|\overline{\mathcal{K}}\mathbf v\|_{\mathcal{M}^{-1}},
 \qquad V_h^{-2}H_0\preceq H\preceq U_h^2H_0.
\end{equation}
No smallness of $q_h$, $\Omega_h$, or $\mathcal{K}-\overline{\mathcal{K}}$ is required.
\end{theorem}
\begin{proof}
The upper bound follows from \eqref{eq:wm-decomposition} and 
Lemmas~\ref{lem:wm-remainder} and \ref{lem:wm-multiplication}:
\[
 \|\mathcal{K}\mathbf v\|_{\mathcal{M}^{-1}}
 \le3\beta\|\overline{\mathcal{K}}\mathbf v\|_{\mathcal{M}^{-1}}
                       +C_R\Omega_h\|\mathbf v\|_{\overline{\mathcal{K}}}
 \le U_h\|\overline{\mathcal{K}}\mathbf v\|_{\mathcal{M}^{-1}}.
\]
 For the other one, let $\mathbf f=\mathcal{K}\mathbf v$, then the decomposition gives
$\overline{\mathcal{K}}\mathbf v=\mathsf D^{-1}(\mathbf f-R \mathbf v)$, and therefore
\[
\|\overline{\mathcal{K}}\mathbf v\|_{\mathcal{M}^{-1}}
\le\frac3\alpha\bigl(\|\mathbf f\|_{\mathcal{M}^{-1}}+C_R\Omega_h\|\mathbf v\|_{\overline{\mathcal{K}}}\bigr)\le\frac3\alpha
       \left(1+\frac{C_Rc_*\Omega_h}{\alpha}\right)\|\mathbf f\|_{\mathcal{M}^{-1}}.
\] 
Here the second inequality uses energy coercivity of $\mathcal{K}$.
Finally, $\|\mathcal K\mathbf v\|_{\mathcal M^{-1}}^2=\mathbf v^\top H\mathbf v$ and $\|\overline{\mathcal K}\mathbf v\|_{\mathcal M^{-1}}^2=\mathbf v^\top H_0\mathbf v$, so squaring gives the matrix inequalities.
\end{proof}

\begin{corollary}
\label{cor:wm-uniform}
Let the side length $L_D$ be fixed, and suppose the local constants
$\nu_-,\nu_+,\alpha,\beta$ can be chosen uniformly over a family
of meshes and stochastic spaces, with $\nu_-,\alpha>0$, and suppose
$\Omega_h\le\Omega$ uniformly. Define
\[
 q=C_Rc_*\Omega,\qquad U=3\beta+q,\qquad
 V=\frac3\alpha(1+q/\alpha).
\]
Then for every member of the family, every $\gamma>0$, and every current
mask,
\begin{equation}\label{eq:wm-uniform-spectrum}
 \operatorname{spec}(P_W^{-1}J)
 \subset[\min\{V^{-2},1\},\max\{U^2,1\}].
\end{equation}
The interval is independent of $h,n_P,N_K,\gamma$, and the mask $D(\mathbf p)$.
\end{corollary}
\begin{proof}
The constants in \eqref{eq:wm-referencebounds} and
\eqref{eq:wm-Rbound} depend only on $d,\nu_-,\nu_+$.
Hence Theorem~\ref{thm:wm-graph} gives uniform constants $U,V$, and
Theorem~\ref{thm:wm-algebra} applies with $c=V^{-2}$ and $C=U^2$.
\end{proof}
\subsection{Perturbation bounds and PCG convergence}
Back to stochastic-major ordering, multiplication in
\eqref{eq:wm-SG} gives the exact identity
\begin{equation}\label{eq:wm-perturb}
\mathcal T=I+\mathcal E,\qquad
 \mathcal E=\sum_{\ell=1}^{N_K-1}G_\ell\otimes B_\ell,
 \qquad B_\ell=M_D^{-1/2}K_\ell K_0^{-1}M_D^{1/2}.
\end{equation}
\begin{proposition}[Small-perturbation enclosure]\label{prop:wm-delta}
If $\|\mathcal E\|_2\le\delta<1$, then
\begin{equation}\label{eq:wm-delta-spectrum}
 \operatorname{spec}(P_W^{-1}J)
 \subset[(1-\delta)^2,(1+\delta)^2],\qquad
 \rm cond_ {\rm sp}(P_W^{-1}J)
 \le\left(\frac{1+\delta}{1-\delta}\right)^2.
\end{equation}
Moreover,
\begin{equation}\label{eq:wm-sum}
 \|\mathcal E\|_2\le\sum_{\ell=1}^{N_K-1}\|G_\ell\|_2\|B_\ell\|_2.
\end{equation}
\end{proposition}
\begin{proof}
The triangle and reverse triangle inequalities give
$(1-\delta)\|\mathbf z\|_2\le\|\mathcal T\mathbf z\|_2\le(1+\delta)\|\mathbf z\|_2$.
Hence $\sigma_-\ge1-\delta$ and $\sigma_+\le1+\delta$, so Lemma~\ref{lem:wm-reduction} gives
$(1-\delta)^2H_0\preceq H\preceq(1+\delta)^2H_0$. Theorem~\ref{thm:wm-algebra} with
$c=(1-\delta)^2\le1\le C=(1+\delta)^2$ then yields \eqref{eq:wm-delta-spectrum}.
For \eqref{eq:wm-sum}, the triangle inequality gives
$\|\mathcal E\|_2\le\sum_{\ell}\|G_\ell\otimes B_\ell\|_2$, and
$\|G_\ell\otimes B_\ell\|_2=\|G_\ell\|_2\|B_\ell\|_2$ because the singular values of a Kronecker product are the products of the singular values of its factors (tensor the two singular value decompositions).
\end{proof}
The following proposition gives computable bounds on $\|B_\ell\|_2$ and $\|\mathcal E\|_2$ from local coefficient values, rather than assuming the desired graph-matrix bound.
Let $K_\ell$ be assembled with scalar weights $a_{\ell,Tq}$, and let
$\mu_{\ell,Tq}=a_{\ell,Tq}/a_{0,Tq}$. Define
\[
 d_{\ell i}=\frac1{|\mathcal Q_i|}
                     \sum_{\mathcal Q_i}\mu_{\ell,Tq},\quad
 b_\ell=\max_{T,q}|\mu_{\ell,Tq}|,\quad
 \omega_{\ell,h}=\max_{i,(T,q)\in\mathcal Q_i}
                     \frac{|\mu_{\ell,Tq}-d_{\ell i}|}{h}.
\]
\begin{proposition}
\label{prop:wm-certificates}
With the matrix assumptions above, signed fluctuation arrays
satisfy
\begin{equation}\label{eq:wm-Bbound}
 \|B_\ell\|_2\le3b_\ell+C_Rc_*\omega_{\ell,h}.
\end{equation}
Alternatively, for the jointly assembled blocks, let
$\tau_h=\max_i\|A_i-I\|_2$. Then
\begin{equation}\label{eq:wm-jointdelta}
 \|\mathcal E\|_2\le3\tau_h+C_Rc_*\Omega_h.
\end{equation}
\end{proposition}
\begin{proof}
Let $D_\ell=\operatorname{diag}(d_{\ell i})$ and
$R_\ell=K_\ell-D_\ell K_0$. The proof of Lemma~\ref{lem:wm-remainder} carries over to the scalar case $n_P=1$, with $A_{Tq}$ and $A_i$ replaced by $\mu_{\ell,Tq}$ and $d_{\ell i}$. Positivity is needed only for
$a_{0,Tq}$, not for $\mu_{\ell,Tq}$. Since $|d_{\ell i}|\le b_\ell$, the argument for \eqref{eq:wm-Dnorms} also applies to $D_\ell$. Hence
\[
 \|R_\ell \mathbf v\|_{M_D^{-1}}
 \le C_R\omega_{\ell,h}\|\mathbf v\|_{K_0},\qquad
 \|D_\ell \mathbf f\|_{M_D^{-1}}\le3b_\ell\|\mathbf f\|_{M_D^{-1}}.
\]
Substituting $\mathbf v=K_0^{-1}M_D^{1/2}\mathbf z$ in the definition of $B_\ell$ gives
\[
 \|B_\ell\|_2
 =\sup_{\mathbf z\ne0}\frac{\|M_D^{-1/2}K_\ell K_0^{-1}M_D^{1/2}\mathbf z\|_2}{\|\mathbf z\|_2}
 =\sup_{\mathbf v\ne0}\frac{\|K_\ell \mathbf v\|_{M_D^{-1}}}{\|K_0\mathbf v\|_{M_D^{-1}}}.
\]
Since $K_\ell=D_\ell K_0+R_\ell$, the two estimates above with $\mathbf f=K_0\mathbf v$, together with the first inequality in \eqref{eq:wm-energy} for $n_P=1$, i.e., $\|\mathbf v\|_{K_0}\le c_*\|K_0\mathbf v\|_{M_D^{-1}}$, give
\[
 \|K_\ell \mathbf v\|_{M_D^{-1}}
 \le3b_\ell\|K_0\mathbf v\|_{M_D^{-1}}+C_R\omega_{\ell,h}\|\mathbf v\|_{K_0}
 \le\bigl(3b_\ell+C_Rc_*\omega_{\ell,h}\bigr)\|K_0\mathbf v\|_{M_D^{-1}}.
\]
Taking the supremum over $\mathbf v\ne0$ proves \eqref{eq:wm-Bbound}.
For \eqref{eq:wm-jointdelta}, we find $\mathcal{K}-\overline{\mathcal{K}}=(\mathsf D-I)\overline{\mathcal{K}}+R$ and apply the same estimates with the joint block norm $\tau_h$. Then, taking the supremum over $\mathbf v$ gives the claim.
%
%
\end{proof}

\begin{corollary}
\label{cor:wm-modefamily}
Suppose, uniformly over the considered meshes and stochastic spaces,
\[
 \sum_{\ell=1}^{N_K-1}\|G_\ell\|_2b_\ell\le\eta<1,
 \qquad
 \sum_{\ell=1}^{N_K-1}\|G_\ell\|_2\omega_{\ell,h}\le\Omega_0<\infty.
\]
Then \eqref{eq:wm-AfromG} satisfies the uniform local hypotheses with
$\alpha=1-\eta$, $\beta=1+\eta$, and~$\Omega=\Omega_0$.
Consequently Corollary~\ref{cor:wm-uniform} applies, whether or not
$3\eta+C_Rc_*\Omega_0<1$.
\end{corollary}
\begin{proof}
By \eqref{eq:wm-AfromG} and the triangle inequality,
$\|A_{Tq}-I\|_2\le\sum_\ell|\mu_{\ell,Tq}|\,\|G_\ell\|_2\le\sum_\ell b_\ell\|G_\ell\|_2\le\eta$.
Since every $G_\ell$ is symmetric, so is $A_{Tq}-I$. Hence its eigenvalues lie in $[-\eta,\eta]$, i.e., $(1-\eta)I\preceq A_{Tq}\preceq(1+\eta)I$.
Linearity of the same nodal averaging gives
$A_i=I+\sum_\ell d_{\ell i}G_\ell$ and
\[
 \|A_{Tq}-A_i\|_2
 \le\sum_\ell\|G_\ell\|_2
                  |\mu_{\ell,Tq}-d_{\ell i}|
 \le h\Omega_0.
\]
\end{proof}

It remains to show how the preconditioner is applied.
\begin{proposition}[Exact current-set inverse]\label{prop:wm-woodbury}
Let $C_D=BB^\top$ with $B\in\mathbb R^{N\times r}$; the columns
need not be independent. 
From the SMW identity
$$(A+UCV)^{-1}=A^{-1}-A^{-1}U\bigl(C^{-1}+VA^{-1}U\bigr)^{-1}VA^{-1}.$$
Here we use $A=H_0, U=B, V=B^\top, C=\gamma^{-1}I_r$,
gives \begin{equation}\label{eq:wm-woodbury}
 P_W^{-1}=(H_0+B\,(\gamma^{-1}I_r)\,B^\top)^{-1}=
 H_0^{-1}-H_0^{-1}B\bigl(\gamma I_r+B^\top H_0^{-1}B\bigr)^{-1}B^\top H_0^{-1}.
\end{equation}
For a binary mask, one may take the columns of $E^\top W^{1/2}$
selected by $D$. Reusing unchanged columns of $B,H_0^{-1}B$, computing added
columns, and discarding removed columns preserves this inverse if the
correction matrix is formed and solved for the current columns.
\end{proposition}
\begin{proof}
The matrix $Z_\gamma=\gamma I_r+B^\top H_0^{-1}B$ is SPD since
$\mathbf q^\top Z_\gamma\mathbf q\ge\gamma\|\mathbf q\|_2^2$ for all $\mathbf q\in\R^r$.
Multiplying the right-hand side of
\eqref{eq:wm-woodbury} by $P_W$ on the left gives
\[
 I+B\left[\gamma^{-1}I_r-Z_\gamma^{-1}
                -\gamma^{-1}B^\top H_0^{-1}BZ_\gamma^{-1}\right](H_0^{-1}B)^\top=I.
\]
For $r=0$, the correction is absent.
Because $H_0$ is fixed, an unchanged column has the same inverse image.
Reordering the columns simply permutes $Z_\gamma$.
Finally, we have
\[
 H_0^{-1}=I_{n_P}\otimes(K_0^{-1}M_DK_0^{-1}),
\]
so its action uses reference stiffness solves and a mass multiplication,
without assembling the double-stiffness matrix.
\end{proof}

\begin{corollary}[PCG consequence]\label{cor:wm-pcg}
Consider one Newton system $J\mathbf x=\mathbf b$ with current-set preconditioner $P_W$.
If \eqref{eq:wm-spectrum} holds, let
$\chi=\theta_+/\theta_-$. With exact arithmetic and exact
preconditioner applications, PCG satisfies
\begin{equation}\label{eq:wm-pcg}
 \|\mathbf x_k-\mathbf x_*\|_J
 \le2\left(\frac{\sqrt{\chi}-1}{\sqrt{\chi}+1}\right)^k
                                      \|\mathbf x_0-\mathbf x_*\|_J.
\end{equation}
Under Proposition~\ref{prop:wm-delta}, the right-hand side is bounded
by $2\delta^k\|\mathbf x_0-\mathbf x_*\|_J$.
\end{corollary}
\begin{proof}
PCG is CG for $P_W^{-1/2}JP_W^{-1/2}$ after the substitution
$\widehat{\mathbf x}=P_W^{1/2}\mathbf x$. Its energy-error norm equals $\|\mathbf x-\mathbf x_*\|_J$.
CG minimizes the SPD quadratic over its affine Krylov space. The
Galerkin orthogonality of that minimizer bounds its relative energy
error by the maximum modulus of any degree-$k$ error polynomial
with value one at zero on the spectral interval. The normalized
Chebyshev polynomial gives the factor in \eqref{eq:wm-pcg}; see,
e.g.,~\cite{elman_finite_2014}.
In the perturbation case,
$\sqrt{\chi}\le(1+\delta)/(1-\delta)$, and substitution gives
the factor $\delta$. 
\end{proof}


\section{Numerical Experiments}\label{sec:numerics}
This section reports two examples.  Example~1 has a manufactured
solution and measures the discretization error together with the solver
cost.  Example~2 has no known solution and measures how the solver cost
depends on the noise amplitude, the stochastic space, the mesh and
$\gamma$.
The stochastic approximation uses total-degree Legendre polynomials in random variables. 
We use continuous bilinear finite elements on uniform Cartesian meshes
with $2^{nc}$ elements in each coordinate direction,
that is, $n_e=2^{nc}$ and $h=2^{1-nc}$ in the notation of Assumption~\ref{ass:wm-local}
, where $nc=4,\ldots,8$.  
All runs start from the zero dual variable and use a normalized KKT
stopping tolerance of $10^{-6}$. The inexact Newton forcing tolerance
is bounded below by $10^{-12}$ and above by $\eta_{\max}=10^{-2}$, and at most
$k_{\max}=100$ Newton steps are allowed in both examples.
{Successful runs additionally satisfy a box-constraint violation not exceeding
$10^{-12}$. In Example~1, they also satisfy a normalized primal--dual gap
not exceeding $10^{-9}$.}
The computations use MATLAB R2023b on the 26-core Intel Xeon Gold 6230R machine.

{We use the notation of \eqref{eq:wm-SG}: $\mathcal M=I_{n_P}\otimes M_D$,
$\overline{\mathcal K}=I_{n_P}\otimes K_0$ and
$H_0=\overline{\mathcal K}\mathcal M^{-1}\overline{\mathcal K}$. The stochastic quadrature integrates $\psi_i\psi_j$
exactly, so $V_\xi^\top W_\xi V_\xi=I_{n_P}$ and
$E^\top WE=I_{n_P}\otimes M_L\approx\mathcal M$.

If every control value is free ($D=I_m$) and $\mathcal K$
is replaced by $\overline{\mathcal K}$, the Jacobian matrix \eqref{eq:J}
becomes approximately $H_0+\gamma^{-1}\mathcal M$.  This sum is not
factorized directly.  Following the matching strategy
of~\cite{Pearson-2011-NAS}, we use
\begin{equation}
\label{P_mean}
\begin{aligned}
 P_{\rm mean}
 &=(\overline{\mathcal K}+\gamma^{-1/2}\mathcal M)\mathcal M^{-1}
  (\overline{\mathcal K}+\gamma^{-1/2}\mathcal M)\\
&=I_{n_P}\otimes\bigl[(K_0+\gamma^{-1/2}M_D)M_D^{-1}(K_0+\gamma^{-1/2}M_D)\bigr].
\end{aligned}
\end{equation}
Expanding~\eqref{P_mean} gives
\[
 P_{\rm mean}=H_0+2\gamma^{-1/2}\overline{\mathcal K}+\gamma^{-1}\mathcal M
 =I_{n_P}\otimes\bigl[K_0M_D^{-1}K_0+2\gamma^{-1/2}K_0+\gamma^{-1}M_D\bigr].
\]
The two outer terms reproduce $H_0+\gamma^{-1}\mathcal M$ exactly.  The
cross term is controlled by
$(\overline{\mathcal K}-\gamma^{-1/2}\mathcal M)\mathcal M^{-1}
(\overline{\mathcal K}-\gamma^{-1/2}\mathcal M)\succeq0$, which gives
\[
 H_0+\gamma^{-1}\mathcal M\preceq P_{\rm mean}
 \preceq2\,(H_0+\gamma^{-1}\mathcal M).
\]
Only the sparse spatial matrix $\widehat K_\gamma=K_0+\gamma^{-1/2}M_D\in\R^{n_h\times n_h}$
is factorized by Cholesky decomposition. Using the identity $(B^\top\otimes A)\operatorname{vec}(X)=\operatorname{vec}(AXB)$, we observe from~\eqref{P_mean} that
\[
 P_{\rm mean}^{-1}\mathbf r=\operatorname{vec}\bigl({\widehat K_\gamma^{-1}}M_D{\widehat K_\gamma^{-1}}R\bigr),
\]
where $\mathbf r=\operatorname{vec}(R)$ with
$R\in\R^{n_h\times n_P}$.


In practice, part
of the control is saturated, and the term $\gamma^{-1}\mathcal M$
overestimates $\gamma^{-1}E^\top WD(\mathbf p_k)E$.  At Newton step $k$, we therefore
replace $\gamma^{-1}$ by $\widehat\vartheta_k/\gamma$, where
\[
 \vartheta_k=\frac{1}{m}\operatorname{Tr}D(\mathbf p_k)\in[0,1],
\]
and $\operatorname{Tr}$ stands for the trace operator.
Next, let
$c_k=\sqrt{\widehat\vartheta_k/\gamma}$ and $\widehat K_k=K_0+c_kM_D$. We then  define an approximation~$ P_k$ of $P_{\text{mean}}$ as follows:
\[
 P_k=(\overline{\mathcal K}+c_k\mathcal M)\mathcal M^{-1}
     (\overline{\mathcal K}+c_k\mathcal M)
 ,\qquad
P_k^{-1}\mathbf r=\operatorname{vec}\bigl(\widehat K_k^{-1}M_D\widehat K_k^{-1}R\bigr).
\]
Observe then that  $P_k=P_{\rm mean}$ when $\widehat\vartheta_k=1$ and $P_k=H_0$ when
$\widehat\vartheta_k=0$.  Note that~$\widehat K_k$ is not recomputed at every
Newton step.  Let $t_k=\vartheta_k/\gamma$, and let $t_{\rm last}$ be its
value at the last refactorization.  We refactorize only at the first Newton
step or when
\[
 |t_k-t_{\rm last}|
 >0.15\max\{|t_k|,|t_{\rm last}|,\sqrt{\varepsilon_{\rm mach}}\}.
\]
In that case $\widehat\vartheta_k=\vartheta_k$; otherwise
$\widehat\vartheta_k$ is the fraction stored at the last refactorization
and the old factor is reused.  The free set enters only through the scalar
$\vartheta_k$, so no information on the location of the free set is needed.


\ifnum\revmode>0
\begin{algorithm}
\caption{Direct semismooth Newton--PCG for the reduced dual
equation~\eqref{eq:F}.}
\label{alg:ssn-pcg}
\begin{algorithmic}
\STATE{Set $\mathbf p_0=0$ and evaluate $F(\mathbf p_0)$ and $D(\mathbf p_0)$.}
\FOR{$k=0,1,2,\dots$}
\STATE{Compute the free fraction $\vartheta_k$. Refactorize $\widehat K_k$ if the
test above is met, otherwise reuse the stored factor. This defines $P_k$.}
\STATE{Solve $J_D\,\delta \mathbf p_k=-F(\mathbf p_k)$ by PCG with preconditioner $P_k$
       to the tolerance $\eta_k$ of~\eqref{eq:ew}.}
\STATE{Choose $\alpha_k$ by Armijo backtracking and set
       $\mathbf p_{k+1}=\mathbf p_k+\alpha_k\,\delta \mathbf p_k$.}
\STATE{Recompute the projection, the active set, the residual $r_{k+1}$ and
       the stopping indicators.}
\IF{$r_{k+1}$ meets the stopping tolerance or $k+1=k_{\max}$}
\STATE{\textbf{break}}
\ENDIF
\ENDFOR
\RETURN $\mathbf y=\mathbf y_d-\Mcal^{-1}\Kcal^\top \mathbf p$ and $\mathbf u=\Proj_{[a,b]}(E\mathbf p/\gamma)$
\end{algorithmic}
\end{algorithm}
\fi

\textbf{Example 1}
We consider a control box-constrained stochastic elliptic optimal control problem
on $D=(-1,1)^2$. 
Let $\xi_1,\xi_2,\xi_3$ be independent random variables
uniformly distributed on $[-1,1]$, and set
\[
\Gamma=[-1,1]^3,
\qquad
\rho(\xi)=\frac{1}{8},
\qquad
\kappa(\xi)=10+\sum_{\ell=1}^{3}\xi_\ell.
\]
Thus $7\leq\kappa(\xi)\leq13$ throughout the parameter domain.
The optimization problem is
\begin{equation}
\label{eq:mms-primal}
\begin{aligned}
\min_{y,u}\quad&
\frac12\mathbb E\!\left[
    \|y-y_d\|_{L^2(D)}^2
\right]
+\frac{\gamma}{2}\mathbb E\!\left[
    \|u\|_{L^2(D)}^2
\right],
\\
\text{subject to}\quad&
-\nabla\cdot\bigl(\kappa(\xi)\nabla y(x,\xi)\bigr)
=f(x,\xi)+u(x,\xi),
\\
&y(x,\xi)=0
    \quad\text{on }\partial D,
\\
&-1\leq u(x,\xi)\leq1
    \quad\text{a.e.\ in }D\times\Gamma,
\end{aligned}
\end{equation}
where $\gamma=10^{-3}$.

To obtain a known optimal solution, define
\begin{equation}
\label{eq:mms-exact}
y^*(x,\xi)=\sin(\pi x_1)\sin(\pi x_2),
p^*(x,\xi)=\frac{\beta}{\kappa(\xi)}s(x),
u^*(x,\xi)
=
\operatorname{Proj}_{[-1,1]}
\left(\frac{p^*(x,\xi)}{\gamma}\right).
\end{equation}
Here $s(x)=\sin(\pi x_1)\sin(\pi x_2)$, so that $y^*=s$, and
$\beta=0.02$.
The corresponding target and source are
\begin{equation}
\label{eq:mms-data}
y_d(x)=(1+2\pi^2\beta)s(x),
\qquad
f(x,\xi)=2\pi^2\kappa(\xi)s(x)-u^*(x,\xi).
\end{equation}
These choices satisfy the state equation and the adjoint relation
\[
-\nabla\cdot\bigl(\kappa(\xi)\nabla p^*(x,\xi)\bigr)
=y_d(x)-y^*(x,\xi),
\]
together with the box-projection optimality condition.
The sign convention is therefore consistent with
$u^*=\operatorname{Proj}_{[-1,1]}(p^*/\gamma)$.

For the accuracy assessment, define
\[
\|v\|_{0,\rho}^2
=
\int_\Gamma\int_D |v(x,\xi)|^2\,dx\,\rho(\xi)\,d\xi,
\]
and report the absolute errors
\[
e_y=\|y_h-y^*\|_{0,\rho},
\qquad
e_p=\|p_h-p^*\|_{0,\rho},
\qquad
e_u=\|u_h-u^*\|_{0,\rho}.
\]
At independent evaluation points, the control is reconstructed as
\[
u_h(x,\xi)
=
\operatorname{Proj}_{[-1,1]}
\left(\frac{p_h(x,\xi)}{\gamma}\right),
\]
The empirical state convergence rate is defined by
\[
\mathrm{EOC}_y({nc})
=
\log_2\left(
\frac{e_y(nc-1)}{e_y(nc)}
\right){,}
\]
where $e_y(nc)$ is the error on the mesh with $2^{nc}$ elements per direction.

\begin{table}[htbp]
\centering
\caption{SSN--PCG performance and absolute $L^2_\rho(\Gamma;L^2(D))$
errors for the manufactured problem. The stochastic degree is fixed at
$P=3$ with $\beta=0.02$. PCG mean/max refers to the iterations per Newton step. Here $\gamma=10^{-3}$, and PCG uses the preconditioner $P_k$ of Section~\ref{sec:numerics}.}
\label{tab:mms-hbase-results-scaled}
\resizebox{\linewidth}{!}{%
\begin{tabular}{clccccccc}
\toprule
$nc$ & $N$ & Newton & PCG total & PCG mean/max & $e_y$ & $e_p$ & $e_u$ & $\mathrm{EOC}_y$ \\
\midrule
4 & 4,500     & 3 & 16 & 5.33/7 & $1.5116\times10^{-2}$ & $4.1864\times10^{-5}$ & $1.6610\times10^{-2}$ & 2.004 \\
5 & 19,220    & 3 & 16 & 5.33/7 & $3.7607\times10^{-3}$ & $1.0550\times10^{-5}$ & $4.1985\times10^{-3}$ & 2.007 \\
6 & 79,380    & 3 & 16 & 5.33/7 & $9.3835\times10^{-4}$ & $2.6878\times10^{-6}$ & $1.0740\times10^{-3}$ & 2.003 \\
7 & 322,580   & 3 & 16 & 5.33/7 & $2.3437\times10^{-4}$ & $8.5077\times10^{-7}$ & $3.3461\times10^{-4}$ & 2.001 \\
8 & 1,300,500 & 3 & 16 & 5.33/7 & $5.8571\times10^{-5}$ & $5.6443\times10^{-7}$ & $2.1656\times10^{-4}$ & 2.001 \\
\bottomrule
\end{tabular}%
}
\end{table}
Observe from~\cref{tab:mms-hbase-results-scaled} that the value of EOC$_y$ for all the mesh sizes is roughly $2$, which suggests that the method achieves approximately second-order convergence in the state. {Moreover, we notice small errors corresponding to the state, the adjoint, and the control on fine meshes.}
On every mesh the solver takes three Newton steps and $16$ PCG
iterations in total, at most seven per step, while $N$ grows by a factor of
$289$ from $4{,}500$ to $1{,}300{,}500$.

{Example~2 is designed for stress test for parameters to verify the preconditioner robustness.}

\textbf{Example 2}
The optimization problem is
\begin{equation}
\label{eq:mms-primal-ex2}
\begin{aligned}
\min_{y,u}\quad&
\frac12\mathbb E\!\left[
    \|y-y_d\|_{L^2(D)}^2
\right]
+\frac{\gamma}{2}\mathbb E\!\left[
    \|u\|_{L^2(D)}^2
\right],
\\
\text{subject to}\quad&
-\nabla\cdot\bigl(\kappa(x,\xi)\nabla y(x,\xi)\bigr)
=u(x,\xi),
\\
&y(x,\xi)=0
    \quad\text{on }\partial D,
\\
&0.3\leq u(x,\xi)\leq1
    \quad\text{a.e.\ in }D\times\Gamma,
\end{aligned}
\end{equation}
and the random coefficient $\kappa(x,\xi)$ is the truncated Karhunen--Lo\`eve (KL) expansion $\kappa(x,\xi)=1+\sqrt{3}\,\sigma
 \sum_{j=1}^{K}\sqrt{\lambda_j}\,g_j(x)\xi_j$ with $\xi_j\sim U[-1,1]$, where $(\lambda_j,g_j)$ are the eigenpairs of the covariance function
 \[
    C_g(x,y)=\sigma^2 \exp \left(-\frac{\left|x_1-y_1\right|}{\ell_1}-\frac{\left|x_2-y_2\right|}{\ell_2}\right) \quad \forall(x,y) \in[-1,1]^2,
\]
where $\sigma$ denotes the standard deviation and ${\ell_1=\ell_2=1}$ are correlation lengths.
The discontinuous target is defined on $D=(-1,1)^2$ as
\[
y_d(x,\xi)=
\begin{cases}
1,&x_1\le 0,\quad x_2\le 0,\\
0,&\text{otherwise}.
\end{cases}
\]
We fix $\gamma=10^{-3}$ and vary the amplitude
$\sigma\in\{0.01,0.02,0.05,0.1,0.2\}$, the pair
$(K,P)\in\{(3,3),(4,4),(6,3)\}$ of the number of retained expansion
terms and the total polynomial degree, and the mesh level $nc=4,\dots,8$.
The total-degree Legendre space has $n_P=\binom{K+P}{P}\in\{20,70,84\}$
modes, and with $n_h=(2^{nc}-1)^2$ interior nodes the SG--FE system has
$N=n_hn_P$ unknowns, from $4{,}500$ for $(K,P,nc)=(3,3,4)$ to
$5{,}462{,}100$ for $(K,P,nc)=(6,3,8)$.
Example~2 tests the solver on a problem without a known solution. 
Since $u\ge0.3$ while $y_d=0$ on three quarters of $D$, the
lower bound is expected to be active on a large part of the domain, and
the active set changes during the iteration.

\begin{table}[!ht]
\centering
\caption{Each entry is Newton / mean PCG (maximum PCG) / solve time in seconds with $\gamma=10^{-3}$. Example~2 with the preconditioner $P_k$. Here $\sigma$ is the amplitude, $K$ is the number of random variables, and $P$ is the  degree of the polynomials.}
\label{tab:factorial-matching-robustness}
\resizebox{\textwidth}{!}{%
\begin{tabular}{c c c c c c c}\toprule
$\sigma$ & $(K,P)$ & $nc=4$ & $nc=5$ & $nc=6$ & $nc=7$ & $nc=8$\\\midrule
0.01 & $(3,3)$ & 8 / 5.50 (9) / 0.162 & 7 / 5.71 (8) / 0.31 & 7 / 5.86 (9) / 1.49 & 7 / 6.43 (10) / 6.78 & 6 / 6.33 (10) / 25.6 \\
  & $(4,4)$ & 9 / 6.22 (11) / 0.444 & 7 / 5.71 (8) / 1.17 & 7 / 6.00 (9) / 5.33 & 7 / 6.29 (10) / 21.7 & 6 / 6.33 (10) / 97.9 \\
  & $(6,3)$ & 8 / 5.38 (8) / 0.845 & 7 / 5.71 (8) / 2.69 & 7 / 6.00 (10) / 10.8 & 6 / 5.83 (9) / 37.5 & 6 / 6.33 (10) / 168 \\
\midrule
0.02 & $(3,3)$ & 8 / 5.75 (9) / 0.175 & 7 / 5.86 (8) / 0.295 & 7 / 6.29 (10) / 1.53 & 7 / 6.86 (10) / 7.11 & 6 / 6.83 (11) / 23.3 \\
  & $(4,4)$ & 8 / 5.62 (9) / 0.38 & 8 / 6.25 (10) / 1.42 & 7 / 6.43 (10) / 5.48 & 6 / 6.33 (10) / 21.7 & 6 / 6.83 (11) / 102 \\
  & $(6,3)$ & 8 / 5.62 (9) / 0.872 & 7 / 5.86 (8) / 2.75 & 7 / 6.43 (10) / 11.1 & 6 / 6.33 (10) / 38.9 & 6 / 6.83 (11) / 155 \\
\midrule
0.05 & $(3,3)$ & 8 / 6.25 (9) / 0.177 & 8 / 6.88 (10) / 0.367 & 7 / 6.86 (10) / 1.61 & 6 / 6.67 (10) / 5.93 & 6 / 7.17 (10) / 27.6 \\
  & $(4,4)$ & 8 / 6.50 (10) / 0.429 & 7 / 6.57 (9) / 1.24 & 7 / 7.00 (11) / 5.67 & 6 / 6.83 (10) / 22.5 & 6 / 7.17 (10) / 104 \\
  & $(6,3)$ & 8 / 6.38 (10) / 0.914 & 7 / 6.43 (10) / 2.83 & 7 / 7.00 (11) / 12.2 & 6 / 6.67 (10) / 39.5 & 6 / 7.17 (10) / 177 \\
\midrule
0.1 & $(3,3)$ & 7 / 6.29 (9) / 0.146 & 7 / 7.29 (11) / 0.343 & 7 / 8.29 (13) / 1.76 & 7 / 8.71 (13) / 8.61 & 6 / 8.50 (13) / 29.6 \\
  & $(4,4)$ & 7 / 6.29 (9) / 0.379 & 7 / 7.29 (11) / 1.32 & 7 / 8.43 (13) / 6.48 & 6 / 8.00 (13) / 24.4 & 6 / 8.67 (13) / 115 \\
  & $(6,3)$ & 8 / 7.00 (12) / 0.984 & 7 / 7.43 (11) / 2.92 & 7 / 8.43 (13) / 13.6 & 6 / 8.00 (13) / 42.9 & 6 / 8.50 (13) / 191 \\
\midrule
0.2 & $(3,3)$ & 7 / 8.57 (12) / 0.18 & 7 / 9.14 (12) / 0.381 & 7 / 10.43 (16) / 2.04 & 6 / 10.00 (14) / 7.41 & 6 / 10.50 (14) / 34.1 \\
  & $(4,4)$ & 8 / 9.25 (15) / 0.506 & 7 / 9.86 (15) / 1.38 & 7 / 10.71 (16) / 7.33 & 6 / 10.17 (14) / 27.7 & 6 / 11.00 (14) / 115 \\
  & $(6,3)$ & 8 / 9.50 (16) / 1.08 & 7 / 9.43 (14) / 4.96 & 6 / 9.67 (14) / 11.6 & 6 / 10.00 (14) / 47.7 & 6 / 10.50 (14) / 212 \\
\bottomrule\end{tabular}%
}
\end{table}

\begin{table}[!ht]
\centering\scriptsize
\caption{Dependence on $\gamma$ of Direct SSN--PCG with the preconditioner $P_k$ at $(K,P,nc,\sigma)=(3,3,4,0.1)$. Entries report Newton steps, total PCG, mean/max PCG, and solve time in seconds.}
\label{tab:organized-gamma-robustness}
\begin{tabular}{c c r r c}
\toprule
$\gamma$  & Newton & PCG & mean/max & time \\
\midrule
 $10^{-1}$  & 4 & 20 & 5.00/7 & 0.074 \\
 $10^{-2}$  & 6 & 36 & 6.00/9 & 0.138 \\
 $10^{-3}$  & 7 & 44 & 6.29/9 & 0.146 \\
 $10^{-4}$  & 16 & 168 & 10.50/16 & 0.364 \\
 $10^{-5}$  & 48 & 937 & 19.52/33 & 2.066 \\
 $10^{-6}$ & 47 & 2,106 & 44.81/127 & 2.980 \\
\bottomrule
\end{tabular}
\end{table}

Table~\ref{tab:factorial-matching-robustness} shows three trends.
First, the Newton count is insensitive to the discretization.  Over all
$75$ runs it lies between $6$ and $9$, and it does not grow under
refinement: it decreases from $7$--$9$ at $nc=4$ to exactly $6$ at $nc=8$
for every $\sigma$ and every $(K,P)$.  This is in line with the
mesh-independent behaviour of semismooth Newton
methods~\cite{Hintermuller-2004-MIS}.

Second, the number of PCG iterations per Newton step is essentially
independent of the mesh and of the stochastic discretization.  For fixed
$\sigma$ and $(K,P)$, the mean count rises by at most about two
iterations from $nc=4$ to $nc=6$ and then stays nearly constant.  For
fixed $\sigma$ and $nc$, the three choices of $(K,P)$ differ by about one
iteration or less, although $n_P$ grows from $20$ to $84$.  

Third, the amplitude $\sigma$ is the only parameter with an effect
on PCG.  At $nc=8$, as $\sigma$ grows from $0.01$ to $0.2$, the mean count
rises from $6.33$ to $10.5$--$11.0$ and the maximum from $10$ to $14$.
This is {consistent with} the theory: the reference $H_0$ uses only the
mean coefficient, and the constants $U$ and $V$ of
Corollary~\ref{cor:wm-uniform} grow with the deviation of the coefficient
from its mean through $\alpha$, $\beta$ and $\Omega$.
Even for
$\sigma=0.2$ the growth is moderate, and the Newton count is unaffected. 
The solution time grows reasonably with the scale of the problem, yet it remains within an acceptable range, even the most difficult one can still be solved in  212 seconds.


Table~\ref{tab:organized-gamma-robustness} separates the two levels of the
method.  Its row $\gamma=10^{-3}$ is the entry $\sigma=0.1$, $(K,P)=(3,3)$,
$nc=4$ of Table~\ref{tab:factorial-matching-robustness}.  
For $\gamma\ge10^{-3}$, the Newton count lies between $4$ and $7$, and
the mean PCG count per Newton step lies between $5.00$ and $6.29$.  Below
$\gamma=10^{-3}$ both counts grow.  At $\gamma=10^{-6}$ the solver takes
$47$ Newton steps with $44.81$ PCG iterations per step on average and $127$
at most, that is, $2{,}106$ PCG iterations in total against $44$ at
$\gamma=10^{-3}$.  This is the regime in which $P_k$ stops being robust.
The growth of the
Newton count is a property of the nonlinear iteration: the projection term
in~\eqref{eq:F} has slope $\gamma^{-1}$, so for small $\gamma$ the reduced
equation approaches a bang-bang system, and more steps are needed before
the fast local convergence of the semismooth Newton method sets in. While, the solve time is acceptable and it grows from $0.074$\,s at
$\gamma=10^{-1}$ to $2.980$\,s at $\gamma=10^{-6}$.

\section{Conclusions}\label{sec:conclusion}
We have studied the optimal control problem of elliptic PDEs with random coefficients subject to a box-constrained control variable via Fenchel duality. In  both the function space and the discrete setting via stochastic Galerkin--finite element methods, the dual formulation yields an unconstrained problem in the adjoint variable and reduces to a single semismooth equation without a duality gap. Our matrix-free semismooth Newton--PCG algorithm solves this equation directly without requiring outer iterations, penalties, or smoothing parameters.
The current-set preconditioner $P_W$ adds the active-set correction to a
mean-based reference operator and is inverted exactly by 
the Sherman--Morrison--Woodbury identity.  For $P_W$ we prove that the
preconditioned spectrum is bounded uniformly in the mesh size, the
stochastic dimension, $\gamma$ and the active set.  
In the numerical experiments, the solver uses the
cheap approximation  $P_k$ of the preconditioner $P_W$. The state solution converges at second order, and the Newton
and PCG counts stay bounded under mesh refinement and enlargement of the
stochastic space.

{Future work includes developing continuation strategies for small regularization parameters, extending the duality framework to lognormal diffusion  coefficients  using weighted spaces, and exploring risk-averse objectives alongside alternative stochastic discretizations.}%
\bibliographystyle{siamplain}
\bibliography{uzawa}


\end{document}